\documentclass[10pt]{article}
\usepackage{amsthm,amsfonts,amsmath}
\usepackage{amscd,amssymb,latexsym,epsfig}
\usepackage[all]{xy}
\title{Noncontractible periodic orbits
       in cotangent bundles
       and Floer homology
       }

\author{Joa Weber\\\\
        ETH-Z\"urich}

      \date{15 April 2004}

\newtheorem{theoremABC}{Theorem}

\newtheorem{theorem}{Theorem}[section]
\newtheorem{corollary}[theorem]{Corollary}

\newtheorem{lemma}[theorem]{Lemma}
\newtheorem{proposition}[theorem]{Proposition}

\newtheorem{definition}[theorem]{Definition}
\newtheorem{remark}[theorem]{Remark}

\newcommand{\C}{{\mathbb{C}}}

\newcommand{\N}{{\mathbb{N}}}

\newcommand{\R}{{\mathbb{R}}}

\newcommand{\Z}{{\mathbb{Z}}}
\newcommand{\Aa}{{\mathcal{A}}}   

\newcommand{\Hh}{{\mathcal{H}}}

\newcommand{\Kk}{{\mathcal{K}}}
\newcommand{\Ll}{{\mathcal{L}}}   
\newcommand{\Mm}{{\mathcal{M}}}   

\newcommand{\Oo}{{\mathcal{O}}}
\newcommand{\Pp}{{\mathcal{P}}}

\newcommand{\Ss}{{\mathcal{S}}}

\newcommand{\Uu}{{\mathcal{U}}}

\newcommand{\Ww}{{\mathcal{W}}}

\newcommand{\im}{{\rm im }}        
\newcommand{\id}{{\rm id}}         

\newcommand{\eps}{{\varepsilon}}

\newcommand{\Cinf}{C^{\infty}}

\newcommand{\inner}[2]{\langle #1, #2\rangle}

\def\NABLA#1{{\mathop{\nabla\kern-.5ex\lower1ex\hbox{$#1$}}}}
\def\Nabla#1{\nabla\kern-.5ex{}_{#1}}
\def\Tabla#1{\Tilde\nabla\kern-.5ex{}_{#1}}
\def\abs#1{\mathopen|#1\mathclose|}   
\def\Abs#1{\left|#1\right|}            

\def\Norm#1{\left\|#1\right\|}

\renewcommand{\Tilde}{\widetilde}

\newcommand{\p}{{\partial}}

\begin{document}

\maketitle


\begin{abstract}
For every nontrivial free homotopy
class $\alpha$ of loops
in any closed connected
Riemannian manifold,
we prove existence
of a noncontractible
1-periodic orbit for every
compactly supported
time-dependent Hamiltonian
on the open unit
cotangent bundle whenever
it is sufficiently large
over the zero section.
The proof shows that the
Biran-Polterovich-Salamon
capacity is finite for every closed connected
Riemannian manifold and every free homotopy
class of loops. This implies
a dense existence theorem for periodic
orbits on level hypersurfaces
and, consequently, a refined version
of the Weinstein conjecture:
Existence of closed characteristics
(one associated
to each nontrivial $\alpha$)
on hypersurfaces in $T^*M$
which are of contact type and
contain the zero section.
\end{abstract}


\section{Introduction and main results}
\label{sec:intro}

Let $M$ be a closed connected smooth
manifold. Let $\pi:T^*M\to M$ be
the cotangent bundle
and $\Omega_{can}=-d\theta$
its canonical symplectic form,
where $\theta$ is the Liouville form.
In addition, we choose a
Riemannian metric on $M$ and
denote the open unit
disc cotangent bundle
by $DT^*M$.
Let $\pi_1(M)$ be the set
of \emph{free}
homotopy classes of loops in $M$
and fix a nontrivial element
$\alpha$.
A Hamiltonian
$H\in C_0^\infty(S^1\times DT^*M)$
gives rise to the Hamiltonian
vector field $X_{H_t}$ on $T^*M$,
which is defined by
$dH_t=\iota(X_{H_t})\Omega_{can}$.
Our aim is to detect
1-periodic orbits of $X_{H_t}$
whose projection to $M$
represents $\alpha$.
This set is denoted by
\begin{equation}\label{eq:1-per-alpha}
     \Pp(H;\alpha)
     :=\{z\in C^\infty(S^1,DT^*M)\mid
     \dot z(t)=X_{H_t}(z(t)), \forall t\in S^1,
     [\pi(z)]=\alpha\}.
\end{equation}
Throughout we identify $S^1$ with $\R/\Z$
and think of $H$
as a smooth, compactly supported
function on $\R\times DT^*M$,
$(t,z)\mapsto H_t(z)$,
satisfying $H_{t+1}=H_t$.
The set of lengths of all periodic
geodesics representing $\alpha$ is
the \emph{marked length spectrum}
$$
     \Lambda_\alpha
     :=\{ length \:(x)
     \mid
     x\in C^\infty(S^1,M),
     \Nabla{t}\p_t x\equiv0,
     [x]=\alpha\}.
$$
The infimum of geodesic lengths
in the class $\alpha$, namely
\begin{equation}\label{eq:ell-alpha}
     \ell_\alpha :=\inf \Lambda_\alpha,
\end{equation}
is indeed realized by a periodic
geodesic according to
Lemma~\ref{le:length-spec} below.

\begin{theoremABC}[Existence of
noncontractible 1-periodic orbit]
\label{thm:exist-orbit}
Let $M$ be a closed connected
smooth Riemannian
manifold and $DT^*M\subset T^*M$
the open unit disc bundle.
Then, for every
nontrivial free homotopy class
$\alpha$ of loops in $M$,
the following is true.
Every Hamiltonian
system $(DT^*M,\Omega_{can},H)$, where
$H\in\Cinf([0,1]\times DT^*M)$
is compactly supported
and satisfies
$$
     \sup_{[0,1]\times M} H
     =:-c\le-\ell_\alpha,
$$
admits a 1-periodic orbit $z$
with $[\pi(z)]=\alpha$
and symplectic action
$\Aa_H(z)\ge c$.
\end{theoremABC}

The existence problem
for periodic orbits has a rich
history (see~\cite{HZ94}, for example).
One of the most
prominent cornerstones
was Floer's solution~\cite{FLOER5}
of the Arnold conjecture
by introducing Floer homology theory.
Periodic orbits detected
in the past were typically
contractible.
It is only very recently that
first steps have been taken
in finding \emph{noncontractible}
ones. Namely, by
Gatien-Lalonde~\cite{GL} and by
Biran-Polterovich-Salamon~\cite{BPS}.
Both approaches require
rather restrictive assumptions
on the manifold, such as flatness.
More precisely,
Theorem~\ref{thm:exist-orbit}
is proved in~\cite{BPS}, if $M$ is
either the flat torus or negatively
curved. We discard these assumptions
\emph{completely}.
Generalizations
of~\cite{GL} are given
in~\cite{Lee}.

Theorem~\ref{thm:exist-orbit}
is sharp in the sense that, firstly,
the inequality cannot be improved
(to see this perturb and smoothen
the Hamiltonian $H(t,x,y)
=\ell_\alpha(-1+\Abs{y})$,
$\Abs{y}\in[0,1]$, appropriately).
Secondly, the zero section
$\Oo_M=M$ in the inequality
cannot be replaced by an
arbitrary smooth
section~\cite[Theorem~C]{BPS}.

\begin{figure}[ht]
\parbox[b]{3cm}{
  \epsfig{figure=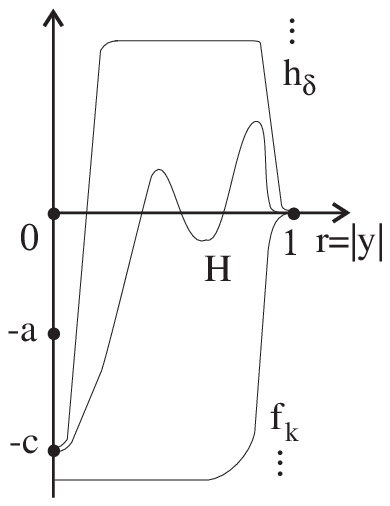,width=\linewidth}
}
\hspace{-.2cm}
\parbox[b]{9.1cm}{
  \begin{equation*}
  \xymatrix{
    &
     \underrightarrow{SH}^{(a,\infty);c;\alpha}_*
     \ar[r]^\simeq
    &
     H_*(\Ll_\alpha M)
    \\
     HF_*^{(a,\infty)}(H;\alpha)
     \ar[ur]^{\iota_H} 
     \ar @{} [r] |{\qquad\circlearrowleft}
    &
     \ar @{} [r] |{\circlearrowleft}
    &
    \\
    &
    \underleftarrow{SH}^{(a,\infty);\alpha}_*
     \ar[ul]^{\pi_H}
     \ar[uu]_T
     \ar[r]^\simeq
    &
     H_*(\Ll^{\frac{a^2}{2}}_\alpha M)
     \ar[uu]_{[\iota]}
  }
  \end{equation*}
}
\caption{Commutative diagram
         in case $a\in(0,c]$.}
\label{fig:fig-idea}
\end{figure}

The proof of
Theorem~\ref{thm:exist-orbit}
is based on combining a main tool
in~\cite{BPS} --
existence of a
natural homomorphism $T$
from symplectic homology of $DT^*M$
to relative symplectic homology
of $(DT^*M,M)$ which
factors through Floer homology
of $H$ (see commuting triangle
in Figure~\ref{fig:fig-idea}) --
and the main result
of~\cite{SW1}
(see also~\cite{Vi98}):
Floer homology of the Hamiltonian
$T^*M\ni(x,y)\mapsto
\frac{1}{2}|y|^2$
is represented by
singular homology of the free
loop space of $M$.
The filtration provided
by the symplectic action functional
is another key ingredient.

A major contribution of the
present text is to compute
symplectic homology
$\underleftarrow{SH}$
of $DT^*M$ and relative
symplectic homology
$\underrightarrow{SH}$
of $(DT^*M,M)$
in terms of the
singular homology of
the free loop space component
$\Ll_\alpha M$
of $M$.
This is the content
of Theorem~\ref{thm:5.1.2}
and it is illustrated by
the rectangular block in
Figure~\ref{fig:fig-idea}.
Here $\iota$
is the natural inclusion
of a sublevel set of
$\Ll_\alpha M$ (see
Section~\ref{sec:FH-convex}).
Note that, in the lower
horizontal isomorphism,
as the action interval $(a,\infty)$
gets bigger homology gets smaller,
because the increasing of the
chain complex admits more
cancellations.
Throughout all homologies are
with coefficients in $\Z_2$.
%
To compute these symplectic
homologies we
construct families of
functions $f_k$ and $h_\delta$,
respectively, as indicated
in Figure~\ref{fig:fig-idea}.
For appropriate $k$ and $\delta$,
they satisfy
$f_k\le H\le h_\delta$.
Together with the composition rule
for the induced
monotone homomorphisms
this leads to
the commutative triangle
in the figure.
Now the map $[\iota]$, hence
$T$, and therefore Floer
homology of $H$
are nonzero
iff $a\in [\ell_\alpha,c]$.
If $c$ happens to be
a regular value of
the symplectic action,
Theorem~\ref{thm:exist-orbit}
follows: Set $a=c$
and use nontriviality
of Floer homology.

Another consequence of
Theorem~\ref{thm:5.1.2} is
finiteness of the
\emph{Biran-Polterovich-Salamon
(BPS) capacity
of $DT^*M$ relative to $M$}
(see Theorem~\ref{thm:3.2.1}
and~(\ref{eq:cBPS}))
\begin{equation*}
\begin{split}
    &c_{BPS}
     (DT^*M,M;\alpha)\\
    &:=\inf\{ c>0\mid
     \Pp(H;\alpha)\not=\emptyset\;\:
     \forall H\in 
     C^\infty_0([0,1]\times DT^*M):
     \sup_{S^1\times M} H\le -c
     \}\\
    &=\ell_\alpha.
\end{split}
\end{equation*}
In \cite{BPS} a
nontrivial free homotopy
class $\alpha$ is called
\emph{symplectically essential}, if
the relative BPS-capacity is finite.
Hence \emph{every}
nontrivial $\alpha$ is
symplectically essential and
this immediately leads to
the following two
\emph{multiplicity results}
(both are proved
in~\cite{BPS} under the assumption
of symplectically essential $\alpha$).

\begin{theoremABC}[Dense existence]
\label{thm:exist-dense}
Let $M$ be a closed connected smooth
manifold and $H:T^*M\to\R$ be a
smooth Hamiltonian which is proper
and bounded from below. Suppose
that the sublevel set $\{ H<c\}$
contains $M$. Then, for every
nontrivial free homotopy class
$\alpha$ of loops in $T^*M$,
there exists a dense
subset $S_\alpha\subset(c,\infty)$
such that the following is true.
For every $s\in S_\alpha$,
the level set $\{ H=s\}$
contains a periodic orbit
$z=(x,y)$ of the Hamiltonian
system $(T^*M,\Omega_{can},H)$
which represents $\alpha$ and
satisfies $\int_0^1
\inner{y(t)}{\dot x(t)} dt>0$.
\end{theoremABC}
\begin{proof}
\cite[Theorem~3.4.1]{BPS}.
\end{proof}

The period
of the orbit in the previous
theorem is not specified.
Note that the theorem
is not true in case $\alpha=0$,
as the example of the flat torus and
Hamiltonian $H(x,y)=\Abs{y}^2/2$ shows.

\begin{theoremABC}
[Closed characteristics]
\label{thm:exist-characteristic}
Let $M$ be a closed connected smooth
manifold and $W\subset T^*M$ be
an open set containing $M$
with compact closure and smooth
convex boundary $Q=\p \overline{W}$.
Let the characteristic line bundle
$\Ll_Q$ be equipped with its canonical
orientation. Then, for every nontrivial
free homotopy class $\alpha$ of loops in $M$,
the characteristic foliation of $Q$
has a closed leaf $z\subset Q$ with
$j_\#[z]=\alpha$. Here $j_\#$ is the map
between free homotopy classes
induced by the composition
$j:Q\subset T^*M\to M$ of inclusion
and projection.
\end{theoremABC}
\begin{proof}
  \cite[Corollary~3.4.2]{BPS}.
\end{proof}

Here convexity of $Q$ means,
by definition, that
there exists a smooth
Liouville vector field $Z$
(i.e. $\Ll_Z\Omega_{can}=\Omega_{can}$),
defined on a neighbourhood $U$ of $Q$
in $\overline{W}$ and pointing
outside $\overline{W}$ along $Q$.
In this case $Q$ is a hypersurface of
\emph{contact type}: $Z$ gives rise
to the contact form $\lambda_Q
:=(\iota(Z)\Omega_{can})\mid_{TQ}$
on $Q$ with $d\lambda_Q
=\Omega_{can}\mid_{TQ}$.
The characteristic line bundle
over $Q$ is given by
$\Ll_Q:=\ker\:(\Omega_{can}\mid_{TQ})$.
The corresponding foliation of $Q$
is called
\emph{characteristic foliation}.
The Reeb vector field $R$ of $\lambda_Q$
is a nonvanishing section of $\Ll_Q$
and induces a
\emph{canonical orientation},
which orients each
leaf of the characteristic foliation.

Theorem~\ref{thm:exist-characteristic}
is related to the celebrated
\emph{Weinstein
conjecture}~\cite{We79}:
Given a symplectic manifold $N$, then
every compact hypersurface $Q\subset N$
of contact type carries a closed
characteristic.
The conjecture was first proved
by Viterbo~\cite{Vi87}
in case $N=\R^{2n}$.
Later Hofer and Viterbo~\cite{HV88}
proved it for cotangent bundles
under the additional hypothesis
$M\subset Q\subset T^*M$.
Hence our \emph{multiplicity} result,
Theorem~\ref{thm:exist-characteristic},
refines their existence theorem.
Viterbo informed us that
the techniques of~\cite{HV88}
should also provide
multiplicities.
In~\cite{Vi97} he proved the
Weinstein conjecture for $T^*M$,
whenever $\pi_1(M)$ is finite.
Further references concerning
Weinstein conjecture and
dense existence
are given in the recent
survey~\cite{Gi03}.
Again, multiplicity results
are obtained only
recently~\cite{GL,BPS,Lee}
under additional assumptions on $M$,
which we dispose of completely.

This paper is organized as follows.
In Section~\ref{sec:floer-homology}
we define action filtered
Floer homology
for a class of Hamiltonians
which are radial and linear
outside a compact subset
of $T^*M$. We study Floer
continuation and give a geometric
criterion to decide if certain
variations of the Hamiltonian
are action-regular (hence
leaving Floer homology invariant).
Theorem~\ref{thm:convex-hamiltonians}
represents Floer homology
for radial convex Hamiltonians
in terms of the loop space of $M$.
In Section~\ref{sec:symp-homology}
we recall the definition
of (relative) symplectic
homology and compute them
in Theorem~\ref{thm:5.1.2}.
In Section~\ref{sec:BPS-capacity}
we compute $c_{BPS}$ and prove
Theorem~\ref{thm:exist-orbit}.

The results of this paper
have been presented at the Workshop
on Symplectic Geometry,
May 12--16 2003, at ETH Z\"urich.

\medskip\noindent{\bf Acknowledgement.}
I wish to thank Dietmar
Salamon for pointing out their
work~\cite{BPS}, as well as for
numerous enlightening
and stimulating discussions.

\section{Floer homology}
\label{sec:floer-homology}

Our aim is to define
Floer homology for
the class of Hamiltonians which,
outside some compact subset of $T^*M$,
are radial and linear
of nonnegative slope.
In Subsection~\ref{sec:FH-convex}
we consider
radial convex Hamiltonians
and relate its Floer homology
to the singular homology
of the free loop space of $M$.

Throughout we use the isomorphism
$T_{(x,y)}T^*M\to T_xM\oplus T_x^*M$
which takes the derivative $\dot z(t)$
of a curve $\R\to T^*M:t\mapsto
z(t)=(x(t),y(t))$,
where $y(t)\in T_{x(t)}^*M$,
to the pair $(\dot x(t),\Nabla{t} y(t))$.
The metric isomorphism
$g:TM\to T^*M$ provides an
almost complex structure $J_g$
and a Riemannian metric on $T^*M$,
both of which are compatible
with $\Omega_{can}$
(see~\cite{SW1}, for example).
For $\alpha\in\pi_1(M)$,
define the space of free loops in $T^*M$
representing $\alpha$ by
$$
     \Ll_\alpha T^*M
     :=\{z=(x,y)\mid \text{$x\in C^\infty(S^1,M)$,
     $[x]=\alpha$, $y(t)\in T_{x(t)}^*M$}\}.
$$
For $H\in C^\infty(S^1\times T^*M)$,
the set $\Pp(H;\alpha)$ of
1-periodic orbits of $H$
representing $\alpha$
corresponds precisely to
the critical points of the
symplectic action functional
$\Aa_H$ restricted to
$\Ll_\alpha T^*M$ and
given by
\begin{equation*}
     \Aa_H(x,y)
     :=\int_0^1\left(
     \left\langle y(t),\dot x(t) \right\rangle
     -H_t(x(t),y(t)) \right) dt.
\end{equation*}
The set of its critical values
is called \emph{action spectrum
of $H$ in the class $\alpha$}
and it is denoted by
$
     Spec(H;\alpha)
     :=\Aa_H(\Pp(H;\alpha))
$.
\begin{remark}[Radial Hamiltonians]
\label{ex:main}\rm
Let $f:\R\to\R$ be smooth and
such that $f(r)=f(-r)$ for
every $r\in\R$. Consider
the Hamiltonian
$H^f:S^1\times T^*M\to \R$, $(t,x,y)
\mapsto f(\abs{y})$.
For simplicity we denote $H^f$
throughout by $f$.
The set of its 1-periodic orbits
can be written as
\begin{equation*}
\begin{split}
     \Pp(f;\alpha)
    &=\big\{(x,y)\in 
     C^\infty(S^1,T^*M)\;\big|\;
     \dot x=\tfrac{f'(\abs{y})}{\abs{y}} y,
     \Nabla{t}y=0, [x]=\alpha\big\} \\
    &=\Pp^+(f;\alpha)\cup
     \Pp^-(f;\alpha)
     \\
     \Pp^\pm(f;\alpha):
    &= \big\{
     z=(x,y)\in C^\infty(S^1,T^*M)
     \;\big|\;
     \text{$x$ periodic geodesic,
     $\ell:=\abs{\dot x}$,} \\
    &\qquad \text{$[x]=\alpha$,
     $y(t)=\pm \tfrac{r_z}{\ell} \dot x(t)$
     where $r_z>0$
     satisfies $f'(r_z)=\pm \ell$}\big\}.
\end{split}
\end{equation*}
Note that $\abs{y(t)}=r_z$ and
$\ell\in\Lambda_\alpha$.
Therefore we call
the elements of $\Lambda_\alpha$
also \emph{critical slopes}.
The symplectic action of
$z^\pm\in\Pp^\pm(f;\alpha)$
with $f'(r_{z^\pm})=\pm\ell$ is
$$
     \Aa_f(z^\pm)
     =f'(r_{z^\pm})r_{z^\pm}-f(r_{z^\pm})
     =\pm\ell r_{z^\pm} -f(r_{z^\pm}).
$$
Given $a\in\R$, there are two methods
to determine if the action
of $z^+\in \Pp^+(f;\alpha)$ is
larger than $a$. \emph{Method 1}:
We have
$\Aa_f(z^+)>a$, if and only if
$f(r_{z^+})$ is located strictly below the line
$r\mapsto -a+r\ell$.
\emph{Method 2} (works for
$z^\pm\in \Pp^\pm(f;\alpha)$):
Draw the tangent
to the graph of $f$ at
the point $r_{z^\pm}$
which satisfies $f'(r_{z^\pm})=\pm\ell$.
This tangent $t(r)=\pm\ell r+b$
intersects the vertical
coordinate axis at
$b=-\Aa_f(z^\pm)$.
Both methods are illustrated
in Figure~\ref{fig:fig-acti}.
\begin{figure}[ht]
  \centering
  \epsfig{figure=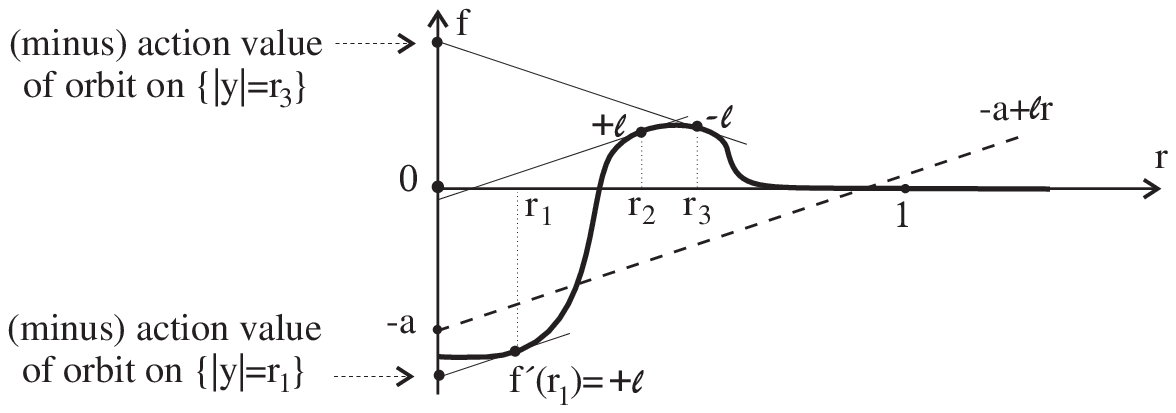}
  \caption{Action of 1-periodic
           orbits on $\{\Abs{y}=r\}$,
           where $f'(r)=\pm\ell
           \in\Lambda_\alpha$.}
  \label{fig:fig-acti}
\end{figure}
\end{remark}

\subsection{Hamiltonians radial
and linear outside a compact set}
\label{subsec:radial-floer-homology}

Consider the class
of smooth time-1-periodic Hamiltonians
on $T^*M$ which are radial and
linear of nonnegative slope
outside some compact subset
of the open disk bundle $D_\rho T^*M$
(of radius $\rho>0$), namely
\begin{equation*}
\begin{split}
     \Kk_\rho
     :=\{ 
     H\in C^\infty(S^1\times T^*M)\mid\;
    &\text{$\exists \eps\in(0,\rho)$
     $\exists \lambda\ge 0$
     $\exists c\in\R$
     such that:} \\
    &\text{$H_t(x,y)=-c+\lambda \Abs{y}$
     whenever $\Abs{y}\ge\rho-\eps$}
     \}.
\end{split}
\end{equation*}
Let $\Kk_\rho$ be equipped
with the $C^\infty$-topology
on the closure of
$S^1\times D_\rho T^*M$.
Fix $\alpha \in \pi_1(M)$
and real numbers
$-\infty\le a<b\le\infty$.
In order to set up Floer
homology for a Hamiltonian $H$ and
restricted to the action window
$(a,b)$ we shall impose
the condition that
$a$ and $b$ are not
in $Spec(H;\alpha)$.
This allows for small
perturbations of $H$.
It is also necessary
to make sure that the
set of relevant periodic orbits
is compact. This is not yet
satisfied:
Consider the case
of $H\in\Kk_\rho$ with
$\lambda \in\Lambda_\alpha$ and
$c\in[a,b]$, which admits 1-periodic
orbits of action $c$ on all sufficiently
large hypersurfaces
$S_\rho T^*M:=\p D_\rho T^*M$.
We simply exclude this situation
by definition, namely
\begin{equation*}
\begin{split}
     \Kk_{\rho;\alpha}^{a,b}
     :=
    &\Bigl\{H\in\Kk_\rho\mid
     \text{
     $\{a,b\}\cap Spec(H;\alpha)
     =\emptyset$ and,
     if $H\equiv-c+\lambda\Abs{y}$
     } \\
    &\text{whenever
     $\Abs{y}\ge\rho$, then
     $\lambda\notin\Lambda_\alpha$
     or $c\notin[a,b]$}
     \Bigr\}.
\end{split}
\end{equation*}

Next we would like to
achieve nondegeneracy
of all elements of the set
$\Pp^{(a,b)}(H;\alpha)$
which, by definition, consists
of all $z\in\Pp(H;\alpha)$
with $\Aa_H(z)\in (a,b)$.
Every 1-periodic
orbit $z$ corresponds to
a fixed point $z_0:=z(0)$
of the time-1-map
$\varphi^H$ of the
Hamiltonian flow
and vice versa.
Recall that $z$ is called
\emph{nondegenerate},
if 1 is not in the spectrum
of $d\varphi^H(z_0):T_{z_0}T^*M
\circlearrowleft$.
It is well known that,
in order to achieve nondegeneracy,
it suffices
to perturb $H$ locally near 
the images of the elements
of $\Pp^{(a,b)}(H;\alpha)$.
A proof in a slightly different
setting may be found
in~\cite[Proof of Theorem~3.6]{JOA1}.
The proof carries over to the
situation at hand and allows
for the following conclusion.
\begin{remark}\label{re:perturb-loop}
\rm
Fix $H\in\Kk_{\rho;\alpha}^{a,b}$
and choose an open neighborhood
$U\subset S^1\times D_{\rho-\eps} T^*M$
of the images of the elements of
$\Pp^{(a,b)}(H;\alpha)$.
Then there exists a neighborhood
$\tilde{\Uu}$ of $0$ in $C_0^\infty(U)$
and a subset $\tilde{\Uu}_{reg}$
of the second category in the
sense of Baire
such that the following is true.
For every $h\in\tilde{\Uu}_{reg}$,
all elements of
$\Pp^{(a,b)}(H+h;\alpha)$
are nondegenerate
and take values in $U$.
Choosing $\tilde{\Uu}_{reg}$
smaller, if necessary,
we may assume without loss
of generality
$H+h\in\Kk_{\rho;\alpha}^{a,b}$,
$\forall h\in\tilde{\Uu}_{reg}$.
\end{remark}
Whereas, for general
$H\in\Kk_{\rho;\alpha}^{a,b}$,
the set $\Pp(H;\alpha)$ may contain 
orbits on all sphere bundles $S_r T^*M$
with $r\in[\rho,\infty)$,
the set $\Pp^{(a,b)}(H;\alpha)$
contains only orbits
taking values in $D_{\rho-\eps} T^*M$.
For all $h\in\tilde{\Uu}_{reg}$,
the latter is also true
in case of $\Pp^{(a,b)}(H+h;\alpha)$,
which is therefore a finite set
by standard arguments.
Now fix $h\in\tilde{\Uu}_{reg}$.
The Conley-Zehnder index
$\mu_{CZ}(z)$ is defined naturally
(see~\cite{JOA1}), for every
orbit $z\in\Pp(H;\alpha)$.
Set $\sigma(z):=0$,
if $x^*TM\to S^1$ is trivial,
and $\sigma(z):=1$ otherwise, and
define the integer
\begin{equation}\label{eq:grading}
     \mu(z):=-\mu_{CZ}(z)+\sigma(z).
\end{equation}
This index provides the grading
for the action filtered
Floer chain groups
\begin{equation}\label{eq:CF}
     CF^{(a,b)}_k(H+h;\alpha)
     :=\underset{\mu(z)=k}
     {\bigoplus_{z\in\Pp^{(a,b)}
     (H+h;\alpha)}} \Z_2 z,
\end{equation}
where $k\in\Z$.
We use the convention that
$CF^{(a,b)}_k(H+h;\alpha):=0$, if
$\Pp^{(a,b)}(H+h;\alpha)=\emptyset$.
Existence of Floer's boundary
operator is based on the following
key proposition (needed at this stage
in its parameter independent form only).

\begin{proposition}[\boldmath$C^0$-bound]
\label{pr:subsolution}
Let $r,c\ge0$
and let $f_s:[r,\infty)\to\R$ 
be a family of smooth functions,
whose dependence on the real parameter
$s$ is of class $C^1$, such that
$\p_s f_s'\ge0$ and
$f_s''\ge -c$ (or $f_s''\le c$),
for all $s\in\R$. Assume that
$H\in C^\infty(\R\times S^1\times TM)$
satisfies
\begin{equation}\label{eq:l1}
     \Abs{y}\ge r \qquad
     \Rightarrow \qquad
     H(s,t,x,y)=f_s(\Abs{y}),
\end{equation}
that the pair
$(u,v)\in C^\infty(\R\times S^1,TM)$
satisfies
\begin{equation}\label{eq:l2}
     \begin{pmatrix}
       \p_su-\Nabla{t}v \\
       \Nabla{s}v+\p_tu
     \end{pmatrix}
     -\nabla H(s,t,u,v) =0,
\end{equation}
and that there exists $T>0$ such that
\begin{equation}\label{eq:l3}
     \Abs{s}\ge T \qquad
     \Rightarrow \qquad
     \Abs{v(s,\cdot)}\le r.
\end{equation}
Then $\Abs{v}\le r$
on $\R\times S^1$.
\end{proposition}

\begin{proof}
Assume by contradiction that there exists
$(s_*,t_*)\in (-T,T)\times S^1$
such that
$$
     r_*
     :=\Abs{v(s_*,t_*)}
     =\max_{[-T,T]\times S^1}\Abs{v}
     >r.
$$
Let $a\in (r,r_*)$ be a regular
value of $\Abs{v}$ and define
the set
$$
     \Omega
     :=\{(s,t)\in\R\times S^1 : \;
     \Abs{v(s,t)}\in [a,r_*]\}
     \subset (-T,T)\times S^1
$$
with smooth boundary
$\partial \Omega=(\Abs{v}^2)^{-1}(a)$.
Denote the connected component
of $\Omega$ containing
$(s_*,t_*)$ by $\Omega_*$.
Assumptions~(\ref{eq:l1}) and~(\ref{eq:l2})
show that
the pair $(u,v)$ restricted to $\Omega_*$
satisfies the equations
$$
     \p_su-\Nabla{t}v=0,\qquad
     \Nabla{s}v+\p_tu
     -\frac{f_s'(\Abs{v})}{\Abs{v}} v=0.
$$
Hence, setting
$L:=\p_s^2+\p_t^2-f_s''(\Abs{v}) \p_s$,
we obtain on $\Omega_*$
$$
     L\tfrac{1}{2}\Abs{v}^2
     =\Abs{\p_su}^2+\Abs{\Nabla{s}v}^2
     +\big(\p_sf_s'\big)(\Abs{v}) \Abs{v}
     \ge 0.
$$
Now the proof, but not the statement,
of the weak maximum 
principle~\cite[Theorem~3.1]{GT77}
carries over to our situation
(we have a bound for $\Abs{s}$ but not
for $\Abs{t}$ when viewed as a real variable):
Consider first the case $f_s''\ge -c$ and
let $\eps\in(0,1]$. Then on $\Omega_*$ we get
$$
     L(\eps e^{-(c+1) s})
     \ge \eps (c+1) e^{-(c+1) s}
     \ge \eps (c+1) e^{-(c+1) T}
     >0.
$$
In the first step we used
$f_s''\ge -c$ and the second step
follows by $s\le T$ on $\Omega_*$.
It follows that the function
$\gamma_\eps:\Omega_*\to\R$,
$(s,t)\mapsto\Abs{v(s,t)}^2
+\eps e^{-(c+1) s}$,
satisfies $L\gamma_\eps>0$
and therefore cannot
have any interior maximum.
(Otherwise, if $(s',t')$ is an interior maximum,
then $\p_s\gamma_\eps(s',t')=0$
and $(\p_s^2+\p_t^2)\gamma_\eps(s',t')\le 0$,
which leads to the contradiction
$(L\gamma_\eps)(s',t') \le 0$).
In other words
$\Norm{\gamma_\eps}_{L^\infty(\Omega_*)}
=\Norm{\gamma_\eps}_{L^\infty(\p\Omega_*)}$
and this proves the third
of the following equations
$$
     r_*^2
     =\Norm{v}_{L^\infty(\Omega_*)}^2
     =\lim_{\eps\to0}
     \Norm{\gamma_\eps}_{L^\infty(\Omega_*)}
     =\lim_{\eps\to0}
     \Norm{\gamma_\eps}_{L^\infty(\p\Omega_*)}
     =\Norm{v}_{L^\infty(\p\Omega_*)}^2
     =a.
$$
To obtain the second equation
we used $s\ge-T$ in $\Omega_*$, which implies
$e^{-(c+1)s}\le e^{(c+1)T}$,
and therefore
$$
     \sup_{\Omega_*} \Abs{v}^2
     \le \sup_{\Omega_*}
     \left(\Abs{v}^2+\eps e^{-(c+1)s}\right)
     \le \left(\sup_{\Omega_*} \Abs{v}^2\right)
     + \eps e^{(c+1)T}
     \stackrel{\eps\to0}{\longrightarrow}
     \sup_{\Omega_*} \Abs{v}^2.
$$
Equation four follows similarly with
$\Omega_*$ replaced by $\p\Omega_*$.
Now $r_*=a$
contradicts the choice of $a$
and this concludes the proof of the lemma
in case $f_s''\ge-c$.
The argument in case $f_s''\le c$ 
is the same up to replacing
$\eps e^{-(c+1) s}$ by $\eps e^{(c+1) s}$.
\end{proof}

Fix $z^\pm\in\Pp^{(a,b)}(H+h;\alpha)$
and consider the moduli space
$\Mm(z^-,z^+;H+h,\alpha)$
of so called
\emph{connecting trajectories}
which consists of
smooth solutions
$w:\R\times S^1\to T^*M$ of
Floer's elliptic partial differential equation
\begin{equation}\label{eq:floer}
     \p_sw +J_g(w)\p_t w -\nabla (H+h)(t,w)=0
\end{equation}
subject to the asymptotic boundary conditions
\begin{equation}\label{eq:limit}
     \lim_{s\to\pm\infty} w(s,t)
     =z^\pm(t), \qquad
     \lim_{s\to\pm\infty} \p_sw(s,t)
     =0,
\end{equation}
uniformly in $t\in S^1$.
Every solution $w$ of~(\ref{eq:floer})
and~(\ref{eq:limit}) satisfies the energy identity
$$
     E(w)
     :=\frac{1}{2}\int_{-\infty}^\infty\int_0^1
     \Abs{w(s,t)}^2 \: dt ds
     =\Aa_{H+h}(z^-)-\Aa_{H+h}(z^+).
$$
More precisely, for solutions $w$
of~(\ref{eq:floer}) which represent
$\alpha$, satisfy~(\ref{eq:l3})
with $r=\rho-\eps$
(i.e. its ends lie in
$D_{\rho-\eps} T^*M$)
and are such that
$\Aa_{H+h}(w(s,\cdot))
\in[a,b]$, for all $s\in\R$,
the following are equivalent:
Existence of the limits~(\ref{eq:limit})
and finite energy of $w$.
Here the key ingredient is
nondegeneracy of all elements of
$\Pp^{(a,b)}(H+h;\alpha)$.
The following observations enable us
to define \emph{Floer's boundary operator}.
\begin{enumerate}
\item[\rm(B1)]
     All relevant Floer connecting trajectories
     stay inside some bounded set
     $D_{\rho-\eps} T^*M$.
     More precisely, for all
     $z^\pm\in \Pp^{(a,b)}(H+h;\alpha)$
     and every $w\in
     \Mm(z^-,z^+;H+h,\alpha)$,
     it holds
     $w(\R\times S^1)\subset
     D_{\rho-\eps} T^*M$:
     All elements of $\Pp^{(a,b)}(H+h;\alpha)$
     lie inside $D_{\rho-\eps} T^*M$,
     for some $\eps>0$.
     Therefore the assumptions
     of Proposition~\ref{pr:subsolution}
     are satisfied for
     $r=\rho-\eps$ and $H+h$ (which
     is radial and linear whenever
     $\Abs{y}\ge \rho-\eps$).
     This proves the claim.
  \item[\rm(B2)]
     The moduli spaces
     $\Mm(z^-,z^+;H+h,\alpha)$
     are compact with respect to
     $C^\infty$-convergence on compact sets.
     This is a consequence of the energy identity,
     exactness of the canonical symplectic form
     $\Omega_{can}$, and~(B1).
  \item[\rm(B3)]
     One can achieve surjectivity
     of the linearized operator
     for equation~(\ref{eq:floer})
     by perturbing $H+h$ in any arbitrarily
     small neighborhood $U$ of the image
     of $w$ in $S^1\times D_{\rho-\eps} T^*M$.
     The perturbation can be chosen
     to vanish up to second order
     along the asymptotic orbits
     $z^\pm$
     (see~\cite[Theorem~5.1~(ii)]{FHS}).
     More precisely, there exists a
     neighborhood $\Uu$ of zero in
     $C_0^\infty(U)$
     and a subset $\Uu_{reg}$
     of the second category
     in the sense of Baire such that
     $\Uu_{reg}\subset
     \tilde{\Uu}_{reg}$ and
     the linearized operator
     for equation~(\ref{eq:floer})
     is surjective, for all
     $w\in \Mm(z^-,z^+;H+h,\alpha)$,
     $z^\pm\in \Pp^{(a,b)}(H+h;\alpha)$
     and $h\in\Uu_{reg}$.
     Hamiltonians satisfying
     this surjectivity
     condition are called
     \emph{regular}.
\end{enumerate}
Changing notation, if necessary,
we may assume without loss of
generality that $H$ is regular.
In this case all moduli spaces
$\Mm(z^-,z^+;H,\alpha)$ are smooth
manifolds of dimension
$\mu_{CZ}(z^+)-\mu_{CZ}(z^-)=\mu(z^-)-\mu(z^+)$
(see~\cite{SALZ} and~\cite{Sa97}) and
they admit a free $\R$-action,
given by shifting the $s$-variable:
$w(s,t) \mapsto w(s+\sigma,t)$.
For $z^-\in\Pp^{(a,b)}(H;\alpha)$
with $\mu(z^-)=k$, define
$$
     \p_k=\p_k(H;\alpha):CF^{(a,b)}_k(H;\alpha)
     \to CF^{(a,b)}_{k-1}(H;\alpha)
$$
by
$$
     \p_k z^- 
     :=\underset{\mu(z^+)=k-1}
     {\sum_{z^+\in\Pp^{(a,b)}(H;\alpha)}}
     \#_2
     \left(\Mm(z^-,z^+;H,\alpha)/\R\right)
     \; z^+.
$$
Here $\#_2$ denotes the number
of elements modulo two.
Define $\p_k z^-:=0$, if the sum
is over the empty set.
Since the index difference is one,
the coefficients in the sum are finite.
This is a consequence of
the compactness property~(B2) of the
1-dimensional moduli space components.
The identity $\p^2=0$ follows
as in Floer's original work~\cite{FLOER5}
using the compactness property (B2) of the
2-dimensional moduli space components.
Action filtered
Floer homology
of the Hamiltonian system
$(T^*M,\Omega_{can},H)$
and with respect to
the free homotopy class $\alpha$
is defined by
$$
     HF_k^{(a,b)}(H;\alpha)
     :=\frac{\ker \p_k}{\im \;\p_{k+1}},\qquad
     k\in\Z.
$$
We suppress the almost
complex structure $J_g$ in
the notation, since
it is fixed
once and for all.

\subsection*{Continuation}
\label{subsec:continuation}

We define filtered
Floer homology
$HF_*^{(a,b)}(H;\alpha)$
for \emph{every}
$H\in\Kk_{\rho;\alpha}^{a,b}$
and show that on
connected components of this set
Floer homology
is independent of $H$.

First observe that the subset
$\Kk_{\rho;\alpha}^{a,b}
\subset\Kk_\rho$ is open:
The second condition in the
definition of
$\Kk_{\rho;\alpha}^{a,b}$,
namely $\lambda\notin\Lambda_\alpha
\vee c\notin[a,b]$ is clearly
open, given
Lemma~\ref{le:length-spec} below.
The first condition
$a,b\notin Spec(H;\alpha)$
is open, because its negation
"$a\in Spec(H;\alpha)$
or $b\in Spec(H;\alpha)$" is closed:
Let $H_\nu$ be a sequence
in $\Kk_\rho$ such that
$a\in Spec(H_\nu;\alpha)$
and which converges
to $H\in\Kk_\rho$ in $C^\infty$.
Hence there
is a sequence
of orbits $z_\nu\in
\Pp(H_\nu;\alpha)$
of action $a$.
This sequence admits
a uniform $C^2$-bound,
which implies existence
of a $C^1$-convergent
subsequence. It follows that
the limit $z$ is a 1-periodic
orbit of $H$ of action $a$
and with $[\pi(z)]=\alpha$.

Now fix $H\in\Kk_{\rho;\alpha}^{a,b}$
and choose a convex neighbourhood
$\Ww$ of $H$ \emph{in}
$\Kk_{\rho;\alpha}^{a,b}$
and regular Hamiltonians
$H^0,H^1\in\Ww$.
Outside $D_{\rho-\eps}T^*M$,
for some $\eps>0$, they
are of the form
$H^k_t(x,y)=-c^k+\lambda^k
\Abs{y}$, for $k=0,1$.
Renaming $H^0$ and $H^1$,
if necessary,
we may assume that
$\lambda^0\le\lambda^1$.
Let $\beta:\R\to[0,1]$
be a smooth nondecreasing
cutoff function which is
zero on $(-\infty,-1]$
and one on $[1,\infty)$.
Consider the smooth homotopy
\begin{equation}\label{eq:stand-cvx}
     \R\to\Ww,\qquad
     s\mapsto H_s:=H^0+\beta(s)
     (H^1-H^0).
\end{equation}
Observe that outside
$D_{\rho-\eps}T^*M$
it is of the form
$$
     f_s(\abs{y})=
     -c^0+\lambda^0\Abs{y}
     +\beta(s)\left(c^0-c^1+\Abs{y}
     (\lambda^1-\lambda^0)\right)
$$
and satisfies $\p_sf_s'(\abs{y})
=\dot\beta(s)(\lambda^1-\lambda^0)
\ge0$.
For $z_k\in\Pp^{(a,b)}(H^k;\alpha)$,
define the moduli space
$\Mm(z_0,z_1;H_s,\alpha)$
of so called
\emph{continuation trajectories}
to be the set of smooth solutions
$w:\R\times S^1\to T^*M$ of
Floer's parameter-dependent equation
\begin{equation}\label{eq:floer-s}
     \p_sw +J_g(w)\p_t w -\nabla H_{s,t}(w)=0
\end{equation}
which satisfy
(uniformly in $t\in S^1$)
the asymptotic boundary conditions
\begin{equation}\label{eq:limit-s}
     \lim_{s\to-\infty} w(s,t)
     =z_0(t), \quad
     \lim_{s\to+\infty} w(s,t)
     =z_1(t), \quad
     \lim_{s\to\pm\infty} \p_sw(s,t)
     =0.
\end{equation}
Moreover, every solution $w$ of~(\ref{eq:floer-s})
and~(\ref{eq:limit-s}) satisfies the energy identity
\begin{equation}\label{eq:s-energy}
     E(w)
      =\Aa_{H^0}(z_0)-\Aa_{H^1}(z_1)
      -\int_{-\infty}^{\infty}\int_0^1
      \left(\p_s H_{s,t}\right)(w(s,t)) \; dt ds.
\end{equation}
The following observations
lead to the definition of
Floer's chain map.
\begin{enumerate}
\item[\rm(C1)]
     All continuation trajectories
     $w\in\Mm(z_0,z_1;H_s,\alpha)$
     stay in $D_{\rho-\eps} T^*M$,
     for all $z_k\in\Pp^{(a,b)}
     (H^k;\alpha)$ and $k=0,1$:
     As observed earlier in~(B1),
     this is true
     for the elements of
     $\Pp^{(a,b)}(H^k;\alpha)$,
     for $k=0,1$.
     Since $\p_sf_s'\ge0$,
     it follows that
     $H_s$ satisfies the assumptions
     of Proposition~\ref{pr:subsolution}
     and this proves the claim.
  \item[\rm(C2)]
     The moduli spaces
     $\Mm(z_0,z_1;H_s,\alpha)$
     are compact with respect to
     $C^\infty$-convergence on compact sets.
     Again this is a consequence of the
     energy identity~(\ref{eq:s-energy}),
     exactness of the canonical symplectic form
     $\Omega_{can}$ and~(C1).
  \item[\rm(C3)]
     A perturbation argument similar
     to~(B3) shows that there is a
     subset (of the second category in
     the sense of Baire)
     of \emph{regular homotopies}
     among all homotopies
     $\R\to\Ww$.
     By definition, this means that
     the linearized operator
     for equation~(\ref{eq:floer-s})
     is surjective,
     for all elements $w$
     of all moduli spaces
     $\Mm(z_0,z_1;H_s,\alpha)$.
\end{enumerate}
By~(C3)
we may assume without 
loss of generality
that $H_s$ is a regular homotopy.
In this case all moduli spaces
$\Mm(z_0,z_1;H_s,\alpha)$ are smooth
manifolds of dimension
$\mu(z_0)-\mu(z_1)$.
By the argument
in~\cite[Section~4.4]{BPS},
we can define
\emph{Floer's continuation map}
$\Phi^{10}:
CF^{(a,b)}_k(H^0;\alpha)\to
CF^{(a,b)}_k(H^1;\alpha)$,
whenever $H^0$ and $H^1$
are \emph{sufficiently close} to
$H\in\Kk_{\rho;\alpha}^{a,b}$.
For $z_0\in\Pp^{(a,b)}(H^0;\alpha)$
with $\mu(z_0)=k$, set
\begin{equation}\label{eq:chain-map}
     \Phi^{10} z_0 
     :=\underset{\mu(z_1)=k}
     {\sum_{z_1\in\Pp^{(a,b)}(H_1;\alpha)}}
     \#_2\left(
     \Mm(z_0,z_1;H_s,\alpha)\right)
     \; z_1.
\end{equation}
That the action
window is preserved
is a consequence
of the energy
identity~(\ref{eq:s-energy}).
The fact that $\Phi^{10}$
is a chain map, i.e.
commutes with
the two boundary operators,
is standard and follows
by investigating the boundary
of the 1-dimensional
components of the moduli spaces
$\Mm(z_0^-,z_1^+;H_s,\alpha)$
(see~\cite[Section~3.4]{Sa97}).
Consequently $\Phi^{10}$
descends to filtered Floer homology.
Moreover, again by standard
arguments, it is independent
of the choice of homotopy
$\R\to\Ww$.

\begin{figure}[ht]
  \centering
  \epsfig{figure=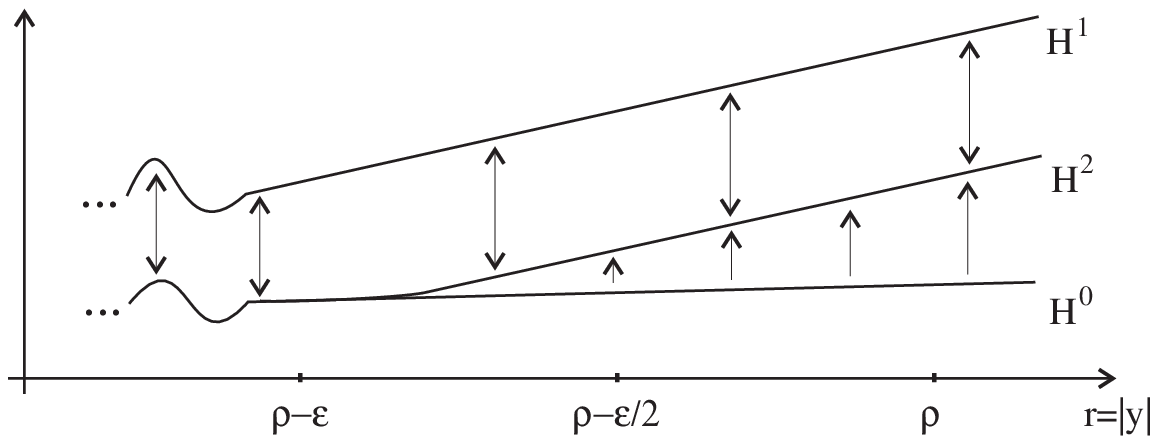}
  \caption{The Hamiltonians
     $H^0,H^1,H^2$ and
     homotopies.}
  \label{fig:fig-hom1}
\end{figure}

In order to prove that
$\Phi^{10}$ induces
an isomorphism on homology,
it is common to consider
the reverse homotopy $H_{-s}$
from $H^1$ to $H^0$
and show that it induces
the inverse map on homology.
Unfortunately, this approach
fails here, because $H_{-s}$
leads to $\p_sf_{-s}'\le0$
and so we can't
justify~$(C1)$ anymore.
We overcome this problem
by introducing
an appropriate intermediate
Hamiltonian $H^2$
(see Figure~\ref{fig:fig-hom1})
and replacing
the homotopy $H^{01}:=H_s$
from $H^0$ to $H^1$ by
$H^{02}_s:=
H^{0}+\beta(s)(H^2-H^0)$
followed by
$H^{21}_s:=
H^2+\beta(s)(H^1-H^2)$.
The corresponding
continuation maps $\Phi^{10}$
and $\Phi^{12}\circ\Phi^{20}$
are chain homotopy equivalent
by standard Floer theory
(see~\cite[Section~3.4]{Sa97})
and therefore the induced
maps on homology are equal.
In other words, the following
diagram commutes
\begin{equation*}
\xymatrix{
    &
     HF^{(a,b)}_*(H^2;\alpha)
     \ar[dr]^{[\Phi^{12}]}
     \ar @{} [d] |{\circlearrowleft}
    &
    \\
     HF^{(a,b)}_*(H^0;\alpha)
     \ar[ur]^{[\Phi^{20}]}
     \ar[rr]^{[\Phi^{10}]}
    &
    &
     HF^{(a,b)}_*(H^1;\alpha).
 }
\end{equation*}
We will prove that $\Phi^{20}=\id$
and\ $\Phi^{12}$ induces
an isomorphism on homology.
Therefore $[\Phi^{10}]$
is an isomorphism.

The Hamiltonian
$H^2$ is defined as follows.
On $D_{\rho-\eps}T^*M$
it coincides with $H^0$
and on the complement
it is radial:
Start at $r=\abs{y}=\rho-\eps$.
Increase $r$ and
make a smooth left turn
and continue, say
outside $D_{\rho-\eps/2}T^*M$,
with the constant slope
$\lambda^1$ of $H^1$
(see Figure~\ref{fig:fig-hom1}).
Proof of $\Phi^{20}=\id$:
If $H^0$ and $H^1$
are sufficiently close,
then the chain groups
associated to $H^0$ and $H^2$
(and $H^1$) are identical.
This follows, since
$\R\setminus\Lambda_\alpha$
is open. Hence all
relevant periodic orbits
lie inside $D_{\rho-\eps}T^*M$. By
Proposition~\ref{pr:subsolution},
this is also true for all
relevant continuation
trajectories~(\ref{eq:floer-s}).
But in $D_{\rho-\eps}T^*M$
the homotopy from $H^0$ to $H^2$
is constant and so they are all
independent of $s$, which
means that the continuation map
$\Phi^{20}$ is the identity.
We show that $\Phi^{12}$
induces an isomorphism
on homology:
The key fact is that the
homotopy $H^{21}$ is outside
$D_{\rho-\eps/2}T^*M$
of the form $f_s(\Abs{y})$
with $\p_sf_s'\equiv 0
\equiv \p_sf_{-s}'$.
Hence, by
Proposition~\ref{pr:subsolution},
all continuation trajectories
corresponding to $H^{21}$
and also those corresponding
to the reverse homotopy
$H^{21}_{-s}$
remain in the bounded
domain $D_{\rho-\eps/2}T^*M$.
Therefore the standard
argument of reversing the homotopy
applies and shows that
$\Phi(H^{21}_{-s})$
induces on homology the inverse
of $\Phi(H^{21}_s)$.
\begin{definition}
\label{def:action-regular}\rm
(1)
  For $H\in\Kk_{\rho;\alpha}^{a,b}$,
  define $HF_*^{(a,b)}(H;\alpha)
  :=HF_*^{(a,b)}(H^0;\alpha)$,
  where $H^0$ is any regular
  Hamiltonian sufficiently close
  to $H$. \\
(2)
  A homotopy $\{H_s\}_{s\in[0,1]}$
  taking values in a connected
  component of
  $\Kk_{\rho;\alpha}^{a,b}$
  is called \emph{action-regular}.
\end{definition}
The observations above
allow to draw the following
conclusions.
Firstly, the
definition in~(1) does not
depend on the choice of $H^0$.
Secondly, given any
action-regular homotopy
$\{H_s\}_{s\in[0,1]}$,
define $s_k:=k/N$ and
choose $N\in\N_0$
sufficiently large to obtain the
composition of isomorphisms
\begin{equation}\label{eq:loc-iso}
     [\Phi^{s_N s_{N-1}}]
     \circ \dots \circ
     [\Phi^{s_1 s_0}]\;:\;
     HF_*^{(a,b)}(H_0;\alpha)
     \to
     HF_*^{(a,b)}(H_1;\alpha).
\end{equation}
Note that it is an open question,
if the resulting isomorphism depends
on the choice of homotopy between
$H_0$ and $H_1$
(see~\cite[Remark~4.4.3]{BPS}).
In the case of a
monotone homotopy
this isomorphism
can be directly defined
in terms of the solutions
of~(\ref{eq:floer-s})
and it is independent of
the choice of the
monotone homotopy.
This will be discussed
in the next section.

\subsection*{Monotone homotopies}
\label{subsec:mon-homotopy}

Given $H^0,H^1\in
\Kk_{\rho;\alpha}^{a,b}$ satisfying
$H^0\le H^1$ pointwise,
there exists a homotopy
$s\mapsto H_s$ from $H^0$
to $H^1$ such that $\p_s H_s\ge 0$.
Such a homotopy is
called \emph{monotone}.
Observe that a monotone
homotopy is outside
$S^1\times D_{\rho-\eps}T^*M$,
for some $\eps>0$,
of the form $f_s(\Abs{y})$
and satisfies necessarily
$\p_sf_s'\ge0$.
By the perturbation
argument~(B3) we may
assume that $H^0$ and $H^1$
are regular.
The energy
identity~(\ref{eq:s-energy})
shows that the map defined
by~(\ref{eq:chain-map})
in terms of the solutions
of~(\ref{eq:floer-s})
is indeed a chain map
from $CF_*^{(a,b)}(H^0;\alpha)$
to $CF_*^{(a,b)}(H^1;\alpha)$
(here the point is
preservation of the action
window).
This shows that a monotone
homotopy induces a
homomorphism, called
\emph{monotone homomorphism},
namely
\begin{equation}\label{eq:mon-hom}
     \sigma^{10}:
     HF_*^{(a,b)}(H^0;\alpha)
     \to
     HF_*^{(a,b)}(H^1;\alpha).
\end{equation}
\begin{lemma}[\cite{FH1,CFH}]\label{le:mon-hom}
  The monotone homomorphism
  does not depend on the
  choice of the monotone homotopy
  used to define it and
  \begin{equation}\label{eq:mon-hom-comp}
    \sigma^{21}\circ\sigma^{10}
    =\sigma^{20},\qquad
    \sigma^{00}=id,
  \end{equation}
  whenever $H^0,H^1,H^2\in
  \Kk_{\rho;\alpha}^{a,b}$
  satisfy $H^0\le H^1\le H^2$.
\end{lemma}
\begin{lemma}[\cite{V}]
\label{le:act-reg-mon-hom}
  The monotone homomorphism
  associated to an action-regular
  monotone homotopy
  is an isomorphism.
\end{lemma}
\begin{proof}
(\cite[Proof of Proposition~4.5.1]{BPS})
The local
isomorphisms~(\ref{eq:loc-iso})
are induced by monotone
homotopies; now
apply~(\ref{eq:mon-hom-comp}).
\end{proof}

Given $a\le b\le c$ and
Hamiltonians $H^0\le H^1\in
\Kk_{\rho;\alpha}^{a,b}\cap
\Kk_{\rho;\alpha}^{a,c}$,
we obtain the following
commuting diagram
(we temporarily
suppress $\alpha$ in the notation),
whose rows are the
short exact sequences
for $H^0$ and for $H^1$
(where $\iota^F$ and $\pi^F$
denote the natural inclusion
and projection, respectively,
and $\alpha$ is fixed)
\begin{equation*}
\begin{split}
\xymatrix{
     0
     \ar[r]
    &
     CF_*^{(a,b)}(H^0)
     \ar[r]^{\iota^F}
     \ar[d]^{\Phi^{10}}
     \ar @{} [dr] |{\circlearrowleft}
    &
     CF_*^{(a,c)}(H^0)
     \ar[r]^{\pi^F}
     \ar[d]^{\Phi^{10}}
     \ar @{} [dr] |{\circlearrowleft}
    &
     CF_*^{(b,c)}(H^0)
     \ar[r]
     \ar[d]^{\Phi^{10}}
    &
     0
    \\
     0
     \ar[r]
    &
     CF_*^{(a,b)}(H^1)
     \ar[r]^{\iota^F}
    &
     CF_*^{(a,c)}(H^1)
     \ar[r]^{\pi^F}
    &
     CF_*^{(b,c)}(H^1)
     \ar[r]
    &
     0
}
\end{split}.
\end{equation*}
The associated long exact sequences
fit into the commutative diagram
\begin{equation}\label{eq:comm-diag1}
\begin{split}
\xymatrix{
     ..
     \ar[r]
    &
     HF_*^{(a,b)}(H^0)
     \ar[r]^{[\iota^F]}
     \ar[d]^{\sigma^{10}}
     \ar @{} [dr] |{\circlearrowleft}
    &
     HF_*^{(a,c)}(H^0)
     \ar[r]^{[\pi^F]}
     \ar[d]^{\sigma^{10}}
     \ar @{} [dr] |{\circlearrowleft}
    &
     HF_*^{(b,c)}(H^0)
     \ar[d]^{\sigma^{10}}
     \ar[r]
    &
     ..
    \\
     ..
     \ar[r]
    &
     HF_*^{(a,b)}(H^1)
     \ar[r]^{[\iota^F]}
    &
     HF_*^{(a,c)}(H^1)
     \ar[r]^{[\pi^F]}
    &
     HF_*^{(b,c)}(H^1)
     \ar[r]
    &
     ..
}
\end{split}.
\end{equation}

Given $a\le b_1\le c$ and
Hamiltonians
$H^0\in\Kk_{\rho;\alpha}^{a,c}$
and $H^1\in
\Kk_{\rho;\alpha}^{a,c}\cap
\Kk_{\rho;\alpha}^{a,b_1}$
with $H^0\le H^1$,
let $H_s$ be a monotone
homotopy between them
such that the induced
continuation map satisfies
$\Phi^{10}
(\Pp^{(a,c)}(H^0;\alpha))
\subset
\Pp^{(-\infty,b_1)}(H^1;\alpha)$.
Then $H_s$ induces a homomorphism
$$
     \Check{\sigma}^{10}:
     HF_*^{(a,c)}(H^0;\alpha)\to
     HF_*^{(a,b_1)}(H^1;\alpha).
$$
If, in addition,
$b_2\in[b_1,c]$,
$\Phi^{10}
(\Pp^{(a,b_2)}(H^0;\alpha))
\subset
\Pp^{(-\infty,b_1)}(H^1;\alpha)$
and
$H^0\in\Kk_{\rho;\alpha}^{a,b_2}$,
then the following diagram
commutes (even on the chain level)
\begin{equation}\label{eq:comm-diag2}
\begin{split}
\xymatrix{
     HF_*^{(a,c)}(H^0;\alpha)
     \ar[r]^{\sigma^{10}}
     \ar[rd]^{\Check{\sigma}^{10}}
    &
     HF_*^{(a,c)}(H^1;\alpha)
    \\
     HF_*^{(a,b_2)}(H^0;\alpha)
     \ar[u]^{[\iota^F]}
     \ar[r]^{\Check{\sigma}^{10}}
    &
     HF_*^{(a,b_1)}(H^1;\alpha)
     \ar[u]_{[\iota^F]}
}
\end{split}.
\end{equation}

\subsection{Computation for radial
convex Hamiltonians}
\label{sec:FH-convex}

Let $\Ll M=C^\infty(S^1,M)$
denote the free loop space of $M$
and $\Ll_\alpha M$ the subspace
whose elements represent
a given free homotopy class $\alpha$
of loops in $M$.
The \emph{classical action
functional} $\Ss_V :\Ll M\to\R$,
for $V\in C^\infty(S^1\times M)$,
is given by
\begin{equation}\label{eq:class-action}
     \Ss_V(x)
     :=\int_0^1\frac{1}{2}
     \Abs{\dot x(t)}^2 -V_t(x(t)) \;dt.
\end{equation}
For $a\in\R$, we define the sublevel set
$
     \Ll_\alpha^a M
     :=\{x\in \Ll_\alpha M\mid
     \Ss_0(x)\le a\}
$
and denote its singular homology
with coefficients in $\Z_2$
by $H_*(\Ll_\alpha^a M)$.
For $b\ge a$, let
$\iota_{ba}:\Ll_\alpha^a M
\hookrightarrow\Ll_\alpha^b M$
denote the natural inclusion.
Set $\iota_a:=\iota_{\infty a}$.
\begin{theorem}[Floer homology
of radial convex Hamiltonians]
\label{thm:convex-hamiltonians}
Fix a free homotopy class $\alpha$.
For every smooth symmetric
function $f$ with $f''\ge0$
the following is true. If
$\lambda$ is a positive
noncritical slope of $f$, i.e.
$\lambda\in\R^+\setminus\Lambda_\alpha$
and $f'(\rho)=\lambda$
for some $\rho>0$,
then there is a natural isomorphism
$$
     \Psi_\lambda^f:
     HF_*^{(-\infty,c_{f,\lambda})}
     (f;\alpha)
     \to
     H_*(\Ll_\alpha^{\lambda^2/2} M),
     \qquad
     c_{f,\lambda}:=\rho f'(\rho)-f(\rho).
$$
If $g$ is another such function,
then the following diagram
commutes
\begin{equation}\label{eq:thm-com1}
\begin{split}
\xymatrix{
     HF_*^{(-\infty,c_{f,\lambda})}
     (f;\alpha)
     \ar[rr]_\simeq
       ^{\Phi_\lambda^{gf}}
     \ar[dr]_{\Psi_\lambda^f}^\simeq
    &
    &
     HF_*^{(-\infty,c_{g,\lambda})}
     (g;\alpha)
     \ar[dl]^{\Psi_\lambda^g}_\simeq
    \\
    &
     H_*(\Ll_\alpha^{\lambda^2/2} M)
    &
}
\end{split}.
\end{equation}
If $\mu\in(0,\lambda]\setminus
\Lambda_\alpha$ is
another slope of $f$,
then the following diagram commutes
\begin{equation}\label{eq:thm-com2}
\begin{split}
\xymatrix{
     HF_*^{(-\infty,c_{f,\mu})}
     (f;\alpha)
     \ar[d]_{\Psi_\mu^f}^\simeq
     \ar[r]^{[\iota^F]}
    &
     HF_*^{(-\infty,c_{f,\lambda})}
     (f;\alpha)
     \ar[d]^{\Psi_\lambda^f}_\simeq
    \\
     H_*(\Ll_\alpha^{{\mu}^2/2} M)
     \ar[r]^{[\iota_{\frac{\lambda^2}{2}
            \frac{\mu^2}{2}}]}
    &
     H_*(\Ll_\alpha^{\lambda^2/2} M)
}
\end{split}.
\end{equation}
\end{theorem}

The proof of the theorem
is based on the main result
in~\cite{SW1}
and the apriori bound provided by
Proposition~\ref{pr:apriori} below.

\begin{proof}
%
Let $f^{(\lambda)}$
be given by following the graph
of $f$ until slope $\lambda$
turns up for the first time,
say at a point $r_\lambda$,
then continue linearly with
slope $\lambda$.
This function is only
$C^1$. After smoothing it out
nearby $r_\lambda$, we may
assume that $f^{(\lambda)}$
is $C^\infty$.
Let $f_0(\Abs{y}):=\Abs{y}^2/2$.
Define
$\Psi_\lambda^f$
to be the composition
of the vertical maps
on the right hand side
of the following
diagram
\begin{equation}\label{eq:diag-1}
\begin{split}
\xymatrix{
     HF_*^{(-\infty,c_{f,\mu})}
     (f+h;\alpha)
     \ar[r]^{[\iota^F]}
     \ar[d]^\equiv
    &
     HF_*^{(-\infty,c_{f,\lambda})}
     (f+h;\alpha)
     \ar[d]_\equiv^{[\id]}
    \\
     HF_*^{(-\infty,\infty)}
     (f^{(\mu)}+\tilde{h};\alpha)
     \ar[r]^{[\iota^F]}_{[\Phi]}
     \ar[d]^\simeq
    &
     HF_*^{(-\infty,\infty)}
     (f^{(\lambda)}+h;\alpha)
     \ar[d]_\simeq^{[\Phi]}
    \\
     HF_*^{(-\infty,\infty)}
     (f_0^{(\mu)}+\tilde{h_0};\alpha)
     \ar[r]^{[\Phi]}_{[\iota^F]}
     \ar[d]^\equiv
    &
     HF_*^{(-\infty,\infty)}
     (f_0^{(\lambda)}+h_0;\alpha)
     \ar[d]_\equiv^{[\id]}
    \\
     HF_*^{(-\infty,\mu^2/2)}
     (f_0+h_0;\alpha)
     \ar[r]^{[\iota^F]}
     \ar[d]^\simeq
    &
     HF_*^{(-\infty,\lambda^2/2)}
     (f_0+h_0;\alpha)
     \ar[d]_\simeq^{[\Phi]}
    \\
     HF_*^{(-\infty,\mu^2/2)}
     (f_0+V;\alpha)
     \ar[r]^{[\iota^F]}
     \ar[d]^\simeq
    &
     HF_*^{(-\infty,\lambda^2/2)}
     (f_0+V;\alpha)
     \ar[d]_\simeq^{\text{
       \scriptsize\cite[Main Theorem]{SW1}}}
    \\
     H_*(\Ll_\alpha^{\mu^2/2} M)
     \ar[r]^{[\iota_{\frac{\lambda^2}{2}
            \frac{\mu^2}{2}}]}
    &
     H_*(\Ll_\alpha^{\lambda^2/2} M)
}
\end{split}.
\end{equation}
Definition of
the vertical maps on the
right hand side
of~(\ref{eq:diag-1}):
By definition,
$$
     HF_*^{(-\infty,c_{f,\lambda})}
     (f;\alpha)
     :=
     HF_*^{(-\infty,c_{f,\lambda})}
     (f+h;\alpha),
$$
where $h\in\Uu_{reg}(f)$ as
in~(B3) is a  small perturbation
of compact support in
$S^1\times D_{r_\lambda-\eps}T^*M$,
for some small $\eps>0$.
As observed repeatedly,
all 1-periodic
orbits of $f+h$
of energy less than
$c_{f,\lambda}$ and all connecting
trajectories (use
Proposition~\ref{pr:subsolution})
lie entirely in
$D_{r_\lambda-\eps}T^*M$.
The same is true for
the Hamiltonian $f^{(\lambda)}+h$
and the action window
$(-\infty,\infty)$.
But in $D_{r_\lambda-\eps}T^*M$
both Hamiltonians coincide,
hence both chain complexes
are identical.

The second vertical map
in~(\ref{eq:diag-1}) is Floer's
continuation map $[\Phi]$ which
is an isomorphism
by the standard \emph{reverse
homotopy argument}
(see Section \emph{Continuation}).
Given $f_0$, let ${f_0}^{(\lambda)}$
and $r_{0,\lambda}$ be
defined as above.
Pick a small perturbation
$h_0\in\Uu_{reg}(f_0)$.
Since the action window is
the whole real line,
we are not restricted
to \emph{monotone} homotopies.
Let the homotopy $H_s$
be a convex combination
of $H_0=f^{(\lambda)}+h$ and
$H_1={f_0}^{(\lambda)}+h_0$
as in~(\ref{eq:stand-cvx}).
Then outside of $D_R T^*M$
with $R:=\max \{r_\lambda,
r_{0,\lambda}\}$
it is of the form $f_s(\Abs{y})$
and satisfies $\p_s f_s'=0$.
Hence $\Phi(H_{-s})$
is defined and induces a map
on homology inverse
to $[\Phi(H_s)]$.
Moreover, any two
convex combinations such as $H_s$
are homotopic through such
convex combinations and therefore
$[\Phi(H_s)]$ is actually independent
of the choice of $H_s$.

The third vertical map 
in~(\ref{eq:diag-1}) is the
identity by the same
argument as the first map.

The fourth vertical map
is again a version Floer's
continuation map. Its
definition and the proof
that it is an isomorphism
require new arguments:
Floer homology as in line five
of diagram~(\ref{eq:diag-1})
is defined in~\cite{SW1}.
The difference is
the perturbation
$V=V(t,x)$ which is a
sufficiently small element
of the set of
\emph{regular potentials}.
This is a subset of the second
category in the sense of Baire of
$C^\infty(S^1\times M)$
(see~\cite[Theorem~1.1]{JOA1}).
For every $(x,y)\in
\Pp^{(-\infty,\lambda^2/2)}
(f_0+V;\alpha)$, it follows
as in~\cite[App.~A]{JOA1} that
\begin{equation*}
\begin{split}
     \Norm{y}_\infty^2
     =\Norm{\dot x}_\infty^2
    &\le \Norm{\dot x}_2^2
     +\Norm{\nabla V}_\infty
     (1+\Norm{\dot x}_2^2) \\
    &\le \lambda^2+2\Norm{V}_\infty
     +\Norm{\nabla V}_\infty
     (1+\lambda^2+2\Norm{V}_\infty).
\end{split}
\end{equation*}
Hence we may assume without loss of
generality, by
choosing $\lambda$
slightly larger in the same
connected component of
$\R\setminus\Lambda_\alpha$
if necessary,
that all relevant 1-periodic orbits
are nondegenerate and
take values in $D_\lambda T^*M$.
The same is true for
the elements of
$\Pp^{(-\infty,\lambda^2/2)}
(f_0+h_0;\alpha)$.
Let $\Phi$ be the continuation
map~(\ref{eq:chain-map}) with
respect to the homotopy
\begin{equation}\label{eq:H-st}
     H_s=H(s,t,x,y)
     :=\tfrac{1}{2}\Abs{y}^2
     +h_0(t,x,y)
     +\beta(s)\left(
     V(t,x)-h_0(t,x,y)\right),
\end{equation}
where $\beta$ is a cutoff function
as in~(\ref{eq:stand-cvx}).
Even though this homotopy is
not monotone, the map
$\Phi$ respects the
action window, whenever
the perturbations $V$ and $h_0$
are sufficiently $C^0$-small
(see~\cite[Section~4.4]{BPS}).
Among conditions~(C1-C3)
only the apriori bound
in~(C1) is nonstandard.
It is provided by
Proposition~\ref{pr:apriori} below,
which says that
there is a uniform
bound $R>\lambda$
for the fibre components
of all relevant continuation
trajectories.
Again the reverse
homotopy argument applies and
shows that $[\Phi]$
is an isomorphism
(see~\cite[Theorem~3.6]{Sa97}).

Existence of
the fifth and last vertical map
in~(\ref{eq:diag-1})
is the main result in~\cite{SW1}.
This proves existence of the
homomorphism $\Phi_\lambda^f$.

Commutativity
of diagram~(\ref{eq:thm-com1})
follows by replacing
the second vertical isomorphism
on the right hand side 
of~(\ref{eq:diag-1})
according to the
following commuting diagram
(invariance of Floer homology;
see~\cite[Theorem~3.7]{Sa97})
\begin{equation*}
\begin{split}
\xymatrix{
     HF_*^{(-\infty,\infty)}
     (f^{(\lambda)}+h;\alpha)
     \ar[d]_{\Phi^{f_0f}}
     \ar[r]^{\Phi^{gf}}
    &
     HF_*^{(-\infty,\infty)}
     (g^{(\lambda)}+h_g;\alpha)
     \ar[dl]^{\Phi^{f_0g}}
    \\
     HF_*^{(-\infty,\infty)}
     (f_0^{(\lambda)}+h_0;\alpha)
    &
 }
\end{split}.
\end{equation*}

To prove the final
claim~(\ref{eq:thm-com2})
we show that
diagram~(\ref{eq:diag-1}) commutes.
The vertical
maps on the left hand side
of~(\ref{eq:diag-1}) are
the analogues of the ones
on the right hand side.
The perturbations on the
left (marked by $\tilde{ }$\;)
are restrictions
of the ones on the right
hand side.
The first five horizontal
maps can be viewed as being
induced by inclusion of
subcomplexes.
Here convexity of the
unperturbed Hamiltonians enters.
Now blocks one, three and four
in diagram~(\ref{eq:diag-1})
commute already on the chain level.
Horizontal maps two and three
can alternatively be interpreted
as being induced by
continuation maps associated
to homotopies.
Note that, due to the fact that
orbits could enter from infinity
during the homotopy,
the induced maps on homology
are in general not isomorphisms.
Continuation maps
satisfy the obvious composition
rule for catenations of homotopies
(see~\cite[Lemmata~3.11 and~3.12]{Sa97})
and this implies commutativity of
block two.
Block five commutes
by~\cite[Main Theorem]{SW1}.
This concludes the proof of
Theorem~\ref{thm:convex-hamiltonians}.
\end{proof}

\begin{proposition}[Apriori bound]
\label{pr:apriori}
Fix constants $c_0,\lambda>0$
and functions
$h\in C^\infty_0(S^1\times D_\lambda T^*M)$
and $V\in C^\infty(S^1\times M)$.
Let $H_{s,t}$ be given by~(\ref{eq:H-st}).
Then there is a constant
$C=C(c_0,h,V)>0$ such that
the following holds.
If $w=(u,v):\R\times S^1\to T^*M$
is a solution of~(\ref{eq:floer-s})
such that
$$
     E(u,v)
     =\Norm{\p_s u}_2^2
     +\Norm{\Nabla{s} v}_2^2
     \le c_0,\qquad
     \sup_{s\in\R} \Aa_{H_{s,t}}
     (u(s,\cdot),v(s,\cdot))
     \le c_0,
$$
then $\Norm{v}_\infty \le C$.
\end{proposition}

\begin{remark}\rm
An apriori bound for
$s$-independent Hamiltonians
was obtained by
Cieliebak~\cite[Theorem~5.4]{Ci94}.
Our proof uses the techniques
developed by Salamon and
the present author in~\cite{SW1}.
There is yet another approach
towards $C^0$-bounds
by Oancea~\cite{Oa03}.
It requires $\p_s H_s\ge0$,
which is not necessarily satisfied
by~(\ref{eq:H-st}).
However,
Proposition~\ref{pr:subsolution}
fits into the
framework~\cite{Oa03}, which
yields some uniform
$C^0$-bound $c$.
In contrast, our proof of
Proposition~\ref{pr:subsolution}
results in the specific value
$c=r$, which is important
for our purposes.
\end{remark}

\begin{proof}
The proof of the apriori estimate
in~\cite[Chapter~4, case $\eps=1$]{SW1}
carries over almost literally.
Just three estimates have to be
verified in the situation at hand.
Identify $T^*M$ and $TM$
via the Riemannian metric.
For $i,j\in\{1,2\}$,
let $\Nabla{i} h_t(x,\xi)\in T_x M$
and
$\Nabla{ij} h_t(x,\xi)\in End(T_x M)$
be determined by
the identities
\begin{equation*}
\begin{split}
     \frac{d}{d\tau} h_t(x,\xi)
    &=:\inner{\Nabla{1} h_t(x,\xi)}{\p_sx}
     +\inner{\Nabla{2} h_t(x,\xi)}
     {\Nabla{s} \xi},\\
     \Nabla{\tau}\Nabla{i} h_t(x,\xi)
    &=:\Nabla{i1} h_t(x,\xi)\p_s x
     +\Nabla{i2} h_t(x,\xi)\Nabla{s} \xi,
\end{split}
\end{equation*}
for every path $\R\to TM:\tau\mapsto
(x(\tau),\xi (\tau))$.
Then equation~(\ref{eq:floer-s}) is equivalent to
\begin{equation}\label{eq:H-st-2}
\begin{split}
     0
    &=\p_su-\Nabla{t} v
     -(1-\sigma)\nabla V_t(u)
     -\sigma\Nabla{1} h_t(u,v),\\
     0
    &=\Nabla{s} v+\p_t u -v
     -\sigma\Nabla{2} h_t(u,v).
\end{split}
\end{equation}
Use these equations to derive
estimate I:
\begin{equation*}
\begin{split}
     (\p_s^2+\p_t^2)\frac{\Abs{v}^2}{2}
    &=\Abs{\Nabla{s} v}^2
     +\Abs{\Nabla{t} v}^2
     +\inner{\Nabla{s}\left(
     v-\p_t u+\sigma\Nabla{2} h_t(u,v)\right)}
     {v}\\
    &\quad+\inner{\Nabla{t}\left(
     \p_s u-(1-\sigma)\nabla V_t(u)
     -\sigma\Nabla{1} h_t(u,v)\right)}{v}\\
    &=\Abs{\Nabla{s} v}^2
     +\Abs{\Nabla{t} v}^2
     +\inner{\Nabla{s} v}{v}\\
    &\quad +\inner{\sigma' \Nabla{2} h_t(u,v)
     +\sigma \Nabla{21} h_t(u,v) \p_su
     +\sigma \Nabla{22} h_t(u,v) \Nabla{s} v}
     {v}\\
    &\quad -(1-\sigma)\inner{
     \Nabla{\p_tu}\nabla V_t(u)
     +\nabla(\p_t V_t)(u)}{v}\\
    &\quad -\sigma\inner{
     \Nabla{11} h_t(u,v) \p_tu
     +\sigma \Nabla{12} h_t(u,v) \Nabla{t} v
     +\Nabla{1} (\p_t h_t)(u,v)}
     {v}\\
%
     &\ge \frac{6}{8}\Abs{\Nabla{s} v}^2
     +\frac{7}{8}\Abs{\Nabla{t} v}^2
     -2\Abs{v}^2\left(\Norm{\Nabla{22} h}_\infty^2
     +\Norm{\Nabla{12} h}_\infty^2\right)\\
    &\quad -\Abs{v}\Bigl(
     2\Norm{\Nabla{2} h}_\infty
     +\Norm{\nabla (\p_tV_t)}_\infty
     +\Norm{\Nabla{1} (\p_th_t)}_\infty\\
    &\quad
     +\Abs{\p_s u} \Norm{\Nabla{21} h}_\infty
     +\Abs{\p_t u} \Norm{\nabla \nabla V}_\infty
     +\Abs{\p_t u} \Norm{\Nabla{11} h}_\infty
     \Bigr)\\
     &\ge \frac{5}{8}\Abs{\Nabla{s} v}^2
     +\frac{6}{8}\Abs{\Nabla{t} v}^2
     -\mu_1\Abs{v}-\mu_2\Abs{v}^2\\
     &\ge \frac{5}{8}\Abs{\Nabla{s} v}^2
     +\frac{6}{8}\Abs{\Nabla{t} v}^2-\frac{1}{2}
     -({\mu_1}^2+2\mu_2)\frac{\Abs{v}^2}{2}
\end{split}
\end{equation*}
where
\begin{equation*}
\begin{split}
     \mu_1
    &=\Norm{\nabla(\p_tV_t)}_\infty
     +\Norm{\Nabla{1}(\p_tV_t)}_\infty
     +\Norm{\Nabla{21} h}_\infty
     \left(\Norm{\nabla V}_\infty
     +\Norm{\Nabla{1} h}_\infty\right)\\
    &\quad +\Norm{\Nabla{2} h}_\infty
     \left(2+\Norm{\nabla\nabla V}_\infty
     +\Norm{\Nabla{11} h}_\infty\right)\\
     \mu_2
    &=2\Norm{\Nabla{22} h}_\infty^2
     +2\Norm{\Nabla{12} h}_\infty^2
     +2\Norm{\Nabla{21} h}_\infty^2\\
    &\quad
     +2\left(
     \Norm{\nabla\nabla V}_\infty
     +\Norm{\Nabla{11} h}_\infty
     \right)^2
     +\Norm{\nabla\nabla V}_\infty
     +\Norm{\Nabla{11} h}_\infty.
\end{split}
\end{equation*}
All $L^\infty$ norms are finite
by compact support of $V$ and $h$.
Here is estimate II:
By assumption, definition of $\Aa_{H_{s,t}}$
and replacing $\p_tu$ according
to~(\ref{eq:H-st-2}),
we obtain for every $s\in\R$
\begin{equation*}
\begin{split}
     c_0
    &\ge\Aa_{H_{s,t}}\left(
     u(s,\cdot),v(s,\cdot)\right)\\
    &=\int_0^1\left(
     \inner{v}{v-\Nabla{s} v+\sigma\Nabla{2} h_t(u,v)}
     -\frac{1}{2}\Abs{v}^2
     -(1-\sigma)V_t(u)
     -\sigma h_t(u,v)
     \right) dt \\
    &\ge
     \int_0^1\left(\frac{1}{4}\Abs{v}^2
     -2\Abs{\Nabla{s} v}^2\right) dt
     -2\Norm{\Nabla{2} h}_\infty^2
     -\Norm{V}_\infty-\Norm{h}_\infty.
\end{split}
\end{equation*}
Estimate III: There is
a constant $a=a(V,h)>0$, such that
$$
     (\p_s^2+\p_t^2-\p_s)\frac{\Abs{v}^2}{2}
     \ge -a\left(\Abs{v}^2+1\right).
$$
This is derived similarly to
estimate I, and $a$ depends
on the $C^2$-norms of $V$
and $h$ similarly to
$\mu_1$ and $\mu_2$.
Given the three estimates,
the proof of
Proposition~\ref{pr:apriori}
proceeds precisely as
in~\cite[Chapter~4, case $\eps=1$]{SW1}.
\end{proof}

\section{(Relative) symplectic
homology of \boldmath$DT^*M$}
\label{sec:symp-homology}

\subsection{Definition}
\label{sec:SH}
To define symplectic homology
we need to introduce direct
and inverse limits. Refering
to~\cite[Section~4.6]{BPS}
for the general case, we shall
discuss these notions only
in the particular situation at hand.
The set $\Hh:=C_0^\infty(S^1\times DT^*M)$
is partially ordered via
$$
     H_0 \preceq H_1 \qquad
     :\Leftrightarrow \qquad
     H_0(t,z)\le H_1(t,z)
     \quad \forall (t,z)\in S^1\times DT^*M.
$$
Fix $\alpha\in\pi_1(M)$, two reals
$-\infty\le a<b\le\infty$
and denote by $\Hh^{a,b}_\alpha$
the set of $H\in\Hh$ such that
$a,b\notin Spec(H;\alpha)$.
The monotone homomorphisms
$\sigma^{H_1H_0}$ of
Section~\ref{sec:floer-homology}
give rise to the \emph{partially ordered
system $(HF,\sigma)$ of
$\Z_2$-vector spaces
over $\Hh^{a,b}_\alpha$}.
By definition, this means
that $HF$ assigns to each
$H\in\Hh^{a,b}_\alpha$ the
$\Z_2$-vector space
${HF_*}^{(a,b)}(H;\alpha)$,
and $\sigma$ assigns to all
elements $H_0\preceq H_1$
of $\Hh^{a,b}_\alpha$ the
monotone homomorphism
$
     \sigma^{H_1H_0}
$
subject to composition
rule~(\ref{eq:mon-hom-comp}).

The partially ordered set
$(\Hh^{a,b}_\alpha,\preceq)$
is \emph{downward directed}:
For all $H_1,H_2\in\Hh^{a,b}_\alpha$
there exists $H_0\in\Hh^{a,b}_\alpha$
such that $H_0\preceq H_1$ and
$H_0\preceq H_2$.
The functor $(HF,\sigma)$
is called an \emph{inverse system
of $\Z_2$-vector spaces
over $\Hh^{a,b}_\alpha$}.
Its \emph{inverse limit}, called
\emph{symplectic homology of $DT^*M$
with respect to $\alpha$ and
the action window $(a,b)$},
is defined by
\begin{equation*}
\begin{split}
    &\underleftarrow{SH}^{(a,b)}_*(DT^*M;\alpha)
     :=\underset{H\in\Hh^{a,b}_\alpha}
     {\underleftarrow{\lim}}
     HF^{(a,b)}_*(H;\alpha) \\
    &:=\left\{\left.\{a_H\}_{H\in\Hh^{a,b}_\alpha}
     \in\prod_{H\in\Hh^{a,b}_\alpha}
     HF^{(a,b)}_*(H;\alpha)\:\right|\:
     H_0\preceq H_1 \Rightarrow
     \sigma^{H_1H_0}(a_{H_0})=a_{H_1} \right\}.
\end{split}
\end{equation*}
For $H\in\Hh^{a,b}_\alpha$, let
\begin{equation}\label{eq:pi}
     \pi_H:
     \underleftarrow{SH}^{(a,b)}_*(DT^*M;\alpha)
     \to
     HF^{(a,b)}_*(H;\alpha)
\end{equation}
be the projection to the
component corresponding to $H$.
It holds
$\pi_{H_1}=\sigma^{H_1H_0}\circ\pi_{H_0}$,
whenever $H_0\preceq H_1$.

To define
\emph{relative} symplectic homology,
fix $c>0$ and consider the subset
$$
     \Hh^{a,b}_{\alpha;c}
     =\Hh^{a,b}_{\alpha;c}(DT^*M,M)
     :=\left\{H\in \Hh^{a,b}_\alpha \mid
     \sup_{S^1\times M} H \le-c\right\}.
$$
This set is upward directed:
For all $H_0,H_1\in\Hh^{a,b}_{\alpha;c}$
there exists $H_2\in\Hh^{a,b}_{\alpha;c}$
such that $H_0\preceq H_2$ and
$H_1\preceq H_2$.
The functor $(HF,\sigma)$
is called a \emph{direct system
of $\Z_2$-vector spaces over
$\Hh^{a,b}_{\alpha;c}$}.
Its \emph{direct limit}, called
\emph{relative symplectic homology of the pair
$(DT^*M,M)$ with respect to
$\alpha$, the action  window
$(a,b)$} and the bound $c$,
is defined to be the quotient
\begin{equation*}
\begin{split}
     \underrightarrow{SH}^{(a,b);c}_*
     (DT^*M,M;\alpha)
    &:=\underset{H\in\Hh^{a,b}_{\alpha;c}}
     {\underrightarrow{\lim}}
     HF^{(a,b)}_*(H;\alpha) \\
    &:=\left\{ (H,a_H)\:\left|\:
     H\in\Hh^{a,b}_{\alpha;c}\:,\:
     a_H\in HF^{(a,b)}_*(H;\alpha)\right.\right\}
     \Big/\sim.
\end{split}
\end{equation*}
Here $(H_0,a_0)\sim(H_1,a_1)$
iff there exists
$H_2\in\Hh^{a,b}_{\alpha;c}$
such that $H_0\preceq H_2$, $H_1\preceq H_2$
and $\sigma^{H_2H_0}(a_0)=\sigma^{H_2H_1}(a_1)$.
This is an equivalence relation, since
$\Hh^{a,b}_{\alpha;c}$
is upward directed.
The direct limit is a
$\Z_2$-vector space
with the operations
$$
     k[H_0,a_0]:=[H_0,ka_0]
     ,\quad
     [H_0,a_0]+[H_1,a_1]
     :=[H_2,\sigma_1{H_2H_0}(a_0)
     +\sigma_{H_2H_1}(a_1)],
$$
for all $k\in\Z_2$ and
$H_2\in\Hh^{a,b}_{\alpha;c}$
such that $H_0\preceq H_2$ and
$H_1\preceq H_2$.
For $H\in\Hh^{a,b}_{\alpha;c}$
define the homomorphism
\begin{equation}\label{eq:iota}
     \iota_H:
     HF^{(a,b)}_*(H;\alpha)
     \to
     \underrightarrow{SH}^{(a,b);c}_*
     (DT^*M,M;\alpha),\qquad
     a_H\mapsto [H,a_H].
\end{equation}
It satisfies
$\iota_{H_0}=\iota_{H_1}\circ
\sigma_{H_1H_0}$,
whenever $H_0\preceq H_1$.

The main feature of symplectic and
relative symplectic homology
for our purposes
is the existence of a unique
homomorphism between them
which factors through Floer homology
(see~\cite[Proposition~4.8.2]{BPS}),
i.e. such that the
diagram
\begin{equation}\label{eq:nat-homo}
\begin{split}
\xymatrix{
     \underleftarrow{SH}^{(a,b)}_*(DT^*M;\alpha)
     \ar[rr]^{T^{(a,b)}_{\alpha;c}}
     \ar[dr]_{\pi_{H'}}
    &
    &
     \underrightarrow{SH}^{(a,b);c}_*
     (DT^*M,M;\alpha)
    \\
    & 
     HF^{(a,b)}_*(H';\alpha)
     \ar[ru]_{\iota_{H'}}
    &
}
\end{split}
\end{equation}
commutes for every
$H'\in\Hh^{a,b}_{\alpha;c}$.
This homomorphism is given by
\begin{equation}\label{eq:def-homo}
     T^{(a,b)}_{\alpha;c}:
     \{a_H\}_{H\in\Hh^{a,b}_\alpha}
     \mapsto
     [H'',a_{H''}],
\end{equation}
where $H''$ is any element of
$\Hh^{a,b}_{\alpha;c}$
(and $T^{(a,b)}_{\alpha;c}$ is
independent of this choice).

\subsection{Computation}
\label{sec:symp-hom-DTM}

\begin{theorem}
[(Relative) symplectic homology of
\boldmath$DT^*M$]
\label{thm:5.1.2}
Let $M$ be a closed smooth Riemannian
manifold and let $\alpha$ be
a free homotopy class of loops in $M$.
Then the following are true.

\noindent
$(i)$ For every $a\in\R\setminus
\Lambda_\alpha$ (with $a>0$, if
$\alpha=0$), there is
a natural isomorphism
$$
     \underleftarrow{SH}^{(a,\infty)}_*
     (DT^*M;\alpha)
     \simeq
     H_*(\Ll^{a^2/2}_\alpha M).
$$
$(ii)$ For all real numbers $a,c>0$, there is
a natural isomorphism
\begin{equation*}
     \underrightarrow{SH}^{(a,\infty);c}_*
     (DT^*M,M;\alpha)
     \simeq
     \begin{cases}
       H_*(\Ll_\alpha M),
      &\text{if $a\in(0,c]$,} \\
       0,
      &\text{if $a>c$.}
     \end{cases}
\end{equation*}
$(iii)$ For every
$a\in(0,c]\setminus \Lambda_\alpha$,
the following
diagram commutes
\begin{equation*}
\begin{split}
\xymatrix{
     \underrightarrow{SH}^{(a,\infty);c}_*
     (DT^*M,M;\alpha)
     \ar[r]^{\quad(ii)}_{\quad\simeq}
    &
     H_*(\Ll_\alpha M)
    \\
     \underleftarrow{SH}^{(a,\infty)}_*
     (DT^*M;\alpha)
     \ar[u]^{T^{(a,\infty)}_{\alpha;c}}
     \ar[r]^{(i)}_\simeq
    &
     H_*(\Ll^{a^2/2}_\alpha M)
     \ar[u]_{[\iota_{a^2/2}]}
}
\end{split}.
\end{equation*}
\end{theorem}

\begin{remark}\rm
1) In part $(i)$ of Theorem~\ref{thm:5.1.2},
assumption $a>0$ in case $\alpha=0$
avoids taking into account the
trivial 1-periodic orbits
near the boundary of $DT^*M$.
Their existence is caused
by the compact support condition
for the Hamiltonians and
their symplectic action is zero.
\\
2) Assumption $a>0$ in part $(ii)$
is to avoid the elements of
$\Pp^-(f;\alpha)$ (see Remark~\ref{ex:main})
whose symplectic action is negative.
\\
3) The homomorphism $T$
in~$(iii)$
is nonzero, if and only if
$a\in[\ell_\alpha,c]$.
\end{remark}

The proof of
Theorem~\ref{thm:5.1.2}
uses the following two lemmata.

\begin{lemma}\label{le:length-spec}
Let $\alpha \in\pi_1(M)$.
The marked length spectrum
$\Lambda_\alpha\subset \R$ is a
closed and nowhere dense subset.
Equivalently, its complement is open and dense.
\end{lemma}

\begin{proof}
The set $\Lambda_\alpha$ is closed:
Let $\{\lambda_j\}_{j\in\N} 
\subset \Lambda_\alpha$ be such that
$\lambda:=\lim_{j\to\infty}\lambda_j$
exists. For each length $\lambda_j$
choose a periodic geodesic
$x_j$ in the class $\alpha$
of length $\lambda_j$.
Since $\Abs{\dot x_j(t)}=\lambda_j$,
$\forall t\in S^1$, there
exists a compact subset
$K\subset TM$ such that
$(q_j,v_j):=(x_j(0),\dot x_j(0))\in K$,
for every $j\in\N$.
Passing to a subsequence, if
necessary, we may assume without loss
of generality that the limit
$$
     (q,v):=
     \lim_{j\to\infty} (q_j,v_j)
$$
exists. Let $x:[0,1]\to M$ be the solution
of the initial value problem
$\Nabla{t}\dot x\equiv0$ and $x(0)=q$,
$\dot x(0)=v$. It follows
that $x$ is periodic:
$$
     (x(0),\dot x(0))
     =\lim_{j\to\infty} (x_j(0),\dot x_j(0))
     =\lim_{j\to\infty} (x_j(1),\dot x_j(1))
     =(x(1),\dot x(1)).
$$
The last step uses the fact that
the solution to a second order
ordinary differential equation
depends continuously on the
initial values.
It remains to show $[x]=\alpha$:
Since $x_j$ converges to $x$
uniformly in $t\in S^1$ and the
injectivity radius of $M$ is positive,
it follows that
$x_j$ is homotopic to $x$,
for all sufficiently large $j$.

The set $\Lambda_\alpha$
does not contain intervals:
It suffices
to prove this property for the set
$
     \Sigma_\alpha
$
of critical values of
the classical action
$\Ss_0:\Ll_\alpha M\to \R$,
since $\Lambda_\alpha=\{\sqrt{2\sigma}\mid \sigma\in
\Sigma_\alpha\}$.
Fix $a\in\R^+$ and consider the
space $L_\alpha^aM$ of piecewise
smooth maps $x:S^1\to M$ such that
$[x]=\alpha$ and $\Ss_0(x)\le a$.
For some sufficiently large
$N=N(a)\in\N$, there is a smooth
finite dimensional manifold
$P_\alpha^N M$
(the set of geodesic polygons;
see~\cite{Bo82}) and
the functional $\Ss_0$ restricted
to $P_\alpha^N M$ is smooth.
The critical points of $\Ss_0$
restricted to the three domains
$$
     P_\alpha^N M
     \subset
     L_\alpha^a M
     \supset
     \Ll_\alpha^a M
$$
all coincide.
Hence
$
     \Sigma_\alpha
     \cap \{\Ss_0(x)\le a\}
$
equals
$
     \{\text{critical values of
     $\Ss_0:P_\alpha^N M\to \R$}\}
$.
The latter set is of
measure zero by the
Theorem of Sard~\cite{Sa42}.

The set $\Lambda_\alpha$
is nowhere dense in $\R$:
By definition of \emph{nowhere dense}
we have to show that every nonempty
open set $U\subset\R$
contains a nonempty
open set $V$ such that
$V\cap\Lambda_\alpha=\emptyset$.
Choose any such $U$ and an
open interval $I\subset\R$
contained in $U$.
Since $\Lambda_\alpha$
does not contain intervals,
we can find an element
$s\in I$ such that
$s\notin \Lambda_\alpha$.
Because $\Lambda_\alpha$ is closed,
it is possible  to choose
an open neighbourhood
$V$ of $s$ in $I$
which is disjoint
from $\Lambda_\alpha$.
\end{proof}

\begin{lemma}\label{le:sing-hom}
For every $\alpha\in\pi_1(M)$,
there is a natural isomorphism
$$
     \underset{a\in\R_0^+}
     {\underrightarrow{\lim}} \;
     H_*(\Ll^a_\alpha M)
     \to H_*(\Ll_\alpha M)
     ,\qquad
     [a,x_a]
     \mapsto [\iota_a](x_a).
$$
\end{lemma}

%
%

\begin{proof}
The direct limit is defined by
$$
     \underset{a\in\R_0^+}
     {\underrightarrow{\lim}} \; 
     H_*(\Ll^a_\alpha M)
     :=\{(a,x_a)\mid
     a\ge 0,\; x_a \in
     H_*(\Ll^a_\alpha M)\} /\sim,
$$
where $(a,x_a)\sim(b,x_b)$
iff there exists $c\ge \max\{ a,b \}$
with
$[\iota_{ca}](x_a)
=[\iota_{cb}](x_b)$.

1) well defined:
Assume $(a,x_a)\sim(b,x_b)$,
i.e. $[\iota_{ca}](x_a)
=[\iota_{cb}](x_b)$
for some $c\ge \max\{ a,b \}$.
Apply $[\iota_c]$ to
both sides to obtain
$
     [\iota_{a}](x_a)
     =[\iota_{b}](x_b)
$.

2) surjective:
Represent $x\in H_*(\Ll_\alpha M)$
by a cycle $f(C)$.
Here $C$ is the total space of
a closed simplicial complex and
$f:C\to\Ll_\alpha M$
is a continuous map.
Let $\Ll_\alpha M$ be equipped
with the $W^{1,2}$-metric, then
$\Ss_0$ given
by~(\ref{eq:class-action})
is continuous
(see~\cite[section~5.4]{JO}).
Since $f(C)$ is compact and
$\Ss_0$ is continuous, we obtain
the real number
$a-1:=\max_{f(C)} \Ss_0$.
Now $f(C) \subset \Ll_\alpha^a M$
shows $x_a:=[f(C)]\in
H_*(\Ll_\alpha^a M)$.

3) injective:
Assuming $[\iota_a](x_a)=0
\in H_*(\Ll_\alpha M)$,
we need to show that
$[a,x_a]$ is the zero element
of the direct limit,
i.e. that there exists
$b\ge a$ such that
$[\iota_{ba}](x_a)
=0\in H_*(\Ll_\alpha^b M)$.
Choose a cycle
$f_a:C\to \Ll_\alpha^a M$
representing $x_a$. By assumption
there is a cycle
$F:W\to \Ll_\alpha M$, where
$W$ is the total space of a
simplicial complex with $\partial W=C$
and $F$ is a continuous map
which coincides with $f_a$ on
$\partial W=C$.
Defining $b-1:=\max_{F(W)}\Ss_0\ge a$
and $F_b:=F : W\to \Ll_\alpha^b M$,
we have
$F_b(p)=F(p)=f_a(p)=(\iota_{ba}
\circ f_a)(p)$, for every
$p\in \partial W=C$.
Hence $[\iota_{ba}](x_a)
=[(\iota_{ba}\circ f_a)(C)]
=[F_b(C)]=0$.
\end{proof}

\begin{proof}[Proof of Theorem~\ref{thm:5.1.2}]
Fix a free homotopy class $\alpha$
of loops in $M$.
We prove Theorem~\ref{thm:5.1.2}
in five steps. Part $(i)$
is equivalent to

\vspace{.15cm}
\noindent
{\sc Step 1.}
{\it $(a)$ For every $a\in\R^+\setminus
\Lambda_\alpha$,
there is a natural isomorphism
$$
     \underleftarrow{SH}
     ^{(a,\infty)}_*(DT^*M;\alpha)
     \simeq H_*(\Ll^{a^2/2}_\alpha M).
$$
$(b)$ If $\alpha\not=0$ and $a\le0$, then
$\underleftarrow{SH}
^{(a,\infty)}_*(DT^*M;\alpha)
=0$.
}

\begin{proof} $(a)$
Consider the 
infimum of all
geodesic lengths
strictly greater than $a$
\begin{equation}\label{eq:ell+}
     \ell_+
     :=\inf\left(\Lambda_\alpha
     \cap(a,\infty)\right).
\end{equation}
Indeed, by Lemma~\ref{le:length-spec},
we have $a<\ell_+\in\Lambda_\alpha$.
We use the convention that
$\inf\emptyset=\infty$.
Consider the sequence of linear functions
$T_k$ (of increasing slopes)
through the points
$(r_{k1},-a)$ and $(r_{k2},0)$, where
\begin{align*}
     r_{k1}
    &=r_{k,1;a,\ell_+}
     :=1-\frac{1}{k}+\frac{3}{16k}
     \left(1-\frac{a}{\ell_+}\right), \\
     r_{k2}
    &=r_{k,2;a,\ell_+}
     :=1-\frac{3}{16k}
     \left(1-\frac{a}{\ell_+}\right),
\end{align*}
and
$$
     T_k(r)=T_{k;a,\ell_+}(r)
     :=ak\frac{r-1+\frac{3}{16k}
     \left(1-\frac{a}{\ell_+}\right)}
     {1-\frac{3}{8}\left(1-\frac{a}{\ell_+}
     \right)}, \qquad k\in \N.
$$
Note that $0<r_{k1}<r_{k2}<1$,
for all $k\in\N$, and that both numbers
converge to $+1$ as $k\to\infty$.
Moreover, the slope $m_1$
of $T_1$ satisfies
\begin{equation}\label{eq:slope}
     a<m_1
     =\frac{a}{1-\frac{3}{8}
     \left(1-\frac{a}{\ell_+}\right)}
     <\ell_+.
\end{equation}
Define the sequence of
piecewise linear functions
(see Figure~\ref{fig:fig-aa-1})
\begin{equation}\label{eq:fam-f}
     \hat{f}_k=\hat{f}_{k;a,\ell_+}(r):=
     \begin{cases}
        -A_k:=T_k(r_{k1}/2),
       &\text{if $r\in[0,r_{k1}/2]$,}\\
        T_k(r),
       &\text{if $r\in[r_{k1}/2,r_{k2}]$,}\\
        0,
       &\text{if $r\ge r_{k2}$}.
     \end{cases}
\end{equation}
\begin{figure}[h]
  \centering
  \epsfig{figure=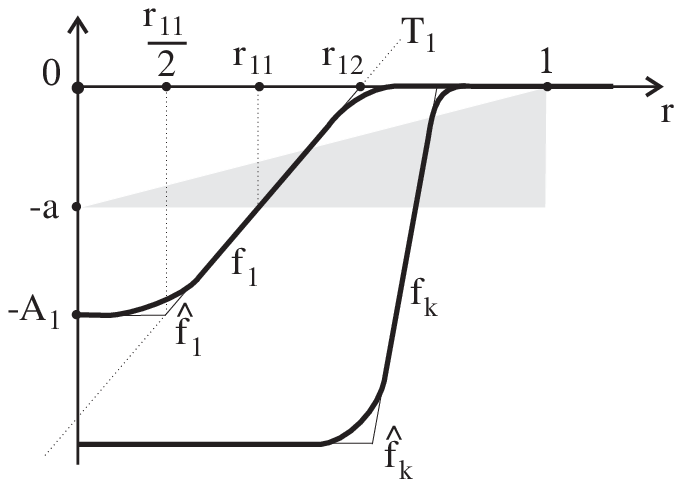}
  \caption{The functions $f_k$.}
  \label{fig:fig-aa-1}
\end{figure}
We obtain a sequence of smooth functions
$f_k=f_{k;a,\ell_+}$ by
smoothing out $\hat{f}_k$
nearby $r_{k1}/2$ and $r_{k2}$
subject to the following conditions.
The new functions $f_k$ coincide
with the old functions $\hat{f}_k$
away from small neighbourhoods
of $r_{k1}/2$ and $r_{k2}$.
In particular, these neighbourhoods
are chosen sufficiently small such that
${\rm graph}\: f_k={\rm graph}\:
\hat{f}_k$ in the
region which lies below the
line $r\mapsto -a+ar$ and above
the line $r\mapsto -a$
(grey shaded in Figure~\ref{fig:fig-aa-1}).
Moreover, the derivatives are required
to satisfy $f_k'\ge0$ everywhere,
$f_k''\ge0$ near $r_{k1}/2$,
$f_k''\le0$ near $r_{k2}$ and
$f_k''=0$ elsewhere
(see Figure~\ref{fig:fig-aa-1}).

The significance of the
functions $f_k$
lies in the fact that
they give rise to natural
isomorphisms
\begin{equation}\label{goal-1}
     \psi_*^{(k)}:
     HF_*^{(a,\infty)}(f_k;\alpha)
     \to
     H_*(\Ll_\alpha^{a^2/2} M),\qquad
     \forall k\in\N.
\end{equation}
Let us first check that, 
$\forall k\in\N$,
indeed $a\notin Spec(f_k;\alpha)$:
The graph of $f_k$
admits a tangent through the
point $(0,-a)$, but
the slope of this tangent
lies in the interval $(a,\ell_+)$
(see Figure~\ref{fig:fig-ff1})
and is therefore noncritical.
\begin{figure}[ht]
  \centering
  \epsfig{figure=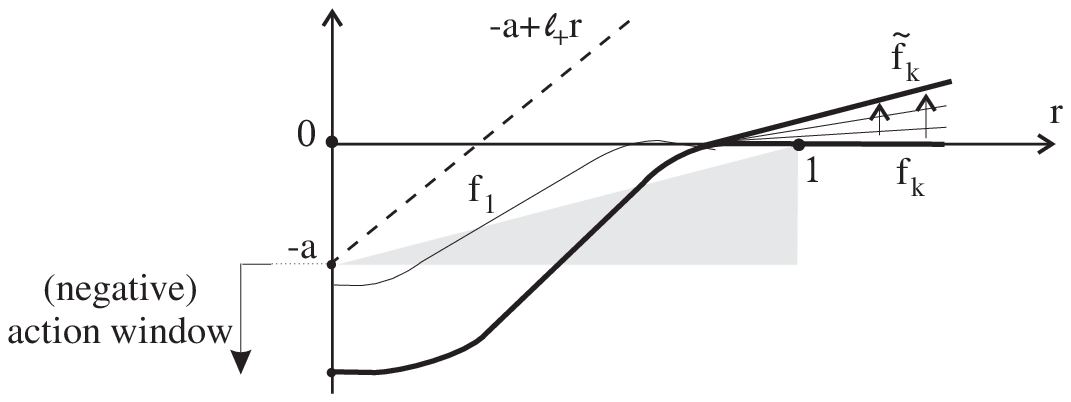}
  \caption{An action-regular monotone
           homotopy between
           $f_k$ and
           $\tilde{f}_k$.}
  \label{fig:fig-ff1}
\end{figure}

Let $\tilde{f}_k$ be the function
obtained by following the
graph of $f_k$ until it
takes on slope $a$
for the second time
(near $r_{k2}$, say at a point $p$),
then continue linearly with
slope $a$. Yet the resulting function
is only of class $C^1$ at $p$.
Local smoothing
near $p$ yields a
smooth function $\tilde{f}_k$
(see Figure~\ref{fig:fig-ff1}).
Figure~\ref{fig:fig-ff1}
also indicates a monotone
action-regular homotopy
which induces the
monotone isomorphism
\begin{equation}\label{eq:psi1}
     \sigma_{\tilde{f}_k f_k}:
     HF_*^{(a,\infty)}(f_k;\alpha)
     \to
     HF_*^{(a,\infty)}(\tilde{f}_k;
     \alpha).
\end{equation}
The homotopy is action-regular,
because all points which
do not remain constant during
the homotopy intersect
the vertical coordinate axis
strictly above the point $-a$.
Therefore all 1-periodic orbits
arising (if any) are of action
strictly less than $a$ (see
Remark~\ref{ex:main}; method 2).

Next consider the function
${f_k}^{(a)}$ obtained by following
the graph of $\tilde{f}_k$
(likewise $f_k$)
until it takes on slope $a$
for the first time (say at a
point $q$), then continue
linearly with slope $a$
(see Figure~\ref{fig:fig-ff2}).
Again smoothing near
$q$ yields ${f_k}^{(a)}$.
\begin{figure}[h]
  \centering
  \epsfig{figure=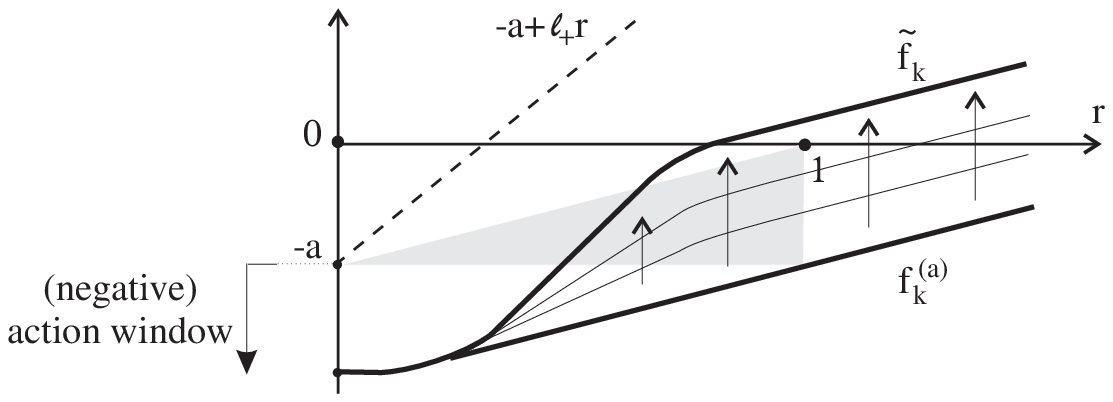}
  \caption{An action-regular monotone
           homotopy between
           ${f_k}^{(a)}$
           and $\tilde{f}_k$.}
  \label{fig:fig-ff2}
\end{figure}
Figure~\ref{fig:fig-ff2}
shows a monotone
homotopy between ${f_k}^{(a)}$
and $\tilde{f}_k$
which is also action-regular:
All tangents to members
of the homotopy which run
through $(0,-a)$
lie above the ray
$\R^+\ni r\mapsto-a+ar$
and strictly below $-a+\ell_+r$.
Hence their slopes
lie in the interval $[a,\ell_+)$
and so are noncritical.
Consequently the homotopy
induces a monotone isomorphism
with inverse
\begin{equation}\label{eq:psi2}
     {\sigma_{\tilde{f}_k f_k^{(a)}}}
     ^{-1}:
     HF_*^{(a,\infty)}(\tilde{f}_k;\alpha)
     \to
     HF_*^{(a,\infty)}({f_k}^{(a)};
     \alpha).
\end{equation}
Since ${f_k}^{(a)}$
does not admit 1-periodic
orbits of action less than $a$,
the map $[\pi^F]$
in~(\ref{eq:comm-diag1}) in
case $(a,b,c)=(-\infty,a,\infty)$
is an isomorphism with inverse
\begin{equation}\label{eq:psi3}
     [\pi^F]^{-1}:
     HF_*^{(a,\infty)}
     ({f_k}^{(a)};\alpha)
     \to
     HF_*^{(-\infty,\infty)}
     ({f_k}^{(a)};\alpha).
\end{equation}
Similarly, there is
the isomorphism
\begin{equation}\label{eq:psi4}
     [\iota^F]^{-1}:
     HF_*^{(-\infty,\infty)}
     ({f_k}^{(a)};\alpha)
     \to
     HF_*^{(-\infty,\tilde{c})}
     ({f_k}^{(a)};\alpha),
\end{equation}
where the constant
$\tilde{c}=\tilde{c}({f_k}^{(a)},a)$
is defined in
Theorem~\ref{thm:convex-hamiltonians}.
Here we need the assumption
that $a$ is not a length.
Theorem~\ref{thm:convex-hamiltonians}
provides the isomorphism
\begin{equation}\label{eq:psi5}
     \Psi^a:
     HF_*^{(-\infty,\tilde{c})}
     ({f_k}^{(a)};\alpha)
     \to
     H_*(\Ll_\alpha^{a^2/2} M).
\end{equation}
Composing~(\ref{eq:psi1}-\ref{eq:psi5})
gives the desired
isomorphisms~(\ref{goal-1}).

To conclude the proof of
part~(a) of Step~1 we show that
$\{f_k\}_{k\in\N}$
is a downward exhausting sequence
(see~\cite[Section~4.7]{BPS})
for the inverse limit
defining symplectic homology.
There are two properties to be checked.
Firstly, given any
$H\in\Hh^{a,\infty}_\alpha$,
there exists $f_k$ such that
$H(t,x,y) \ge f_k(\Abs{y})$,
for every $(t,x,y)\in S^1\times DT^*M$.
This can be achieved by choosing
$k\in\N$ sufficiently large.
Secondly, for every $k\in\N$,
the map
$$
     \sigma_{f_k f_{k+1}}:
     HF_*^{(a,\infty)}(f_{k+1};\alpha)
     \to HF_*^{(a,\infty)}(f_k;\alpha)
$$
induced by a monotone homotopy
between $f_{k+1}$ and $f_k$
is an isomorphism.
To see this
think of $k$ as a continuous
parameter in $[1,\infty)$
instead of a positive integer.
Consider the family
$\{f_r\}_{r\in[1,\infty)}$
defined analogous to
$\{f_k\}_{k\in\N}$.
Then the homotopy
$
     \{f_{k-s}\}
     _{s \in[-1,0]}
$
from $f_{k+1}$ to $f_k$
is clearly monotone and
in addition action-regular:
Similar arguments as above
show that all tangents to members
of the homotopy which run
through $(0,-a)$ are located
between the lines $-a+ar$ and
$-a+\ell_+r$. Hence their
slope is in the interval
$[a,\ell_+)$ which does not
contain critical slopes.
Consequently the boundary $+a$
of the action window remains
disjoint from the action spectrum
throughout the homotopy.
This proves property two
and therefore that the
sequence is exhausting.
Hence the homomorphism~(\ref{eq:pi}),
given in the situation at hand by
\begin{equation}\label{eq:proj-symp}
     \pi_{f_k} : 
     \underleftarrow{SH}^{(a,\infty)}_*
     (DT^*M;\alpha)
     \to
     HF_*^{(a,\infty)}(f_k;\alpha),
\end{equation}
is an isomorphism
by~\cite[Lemma~4.7.1 (ii)]{BPS}.
Then~(\ref{goal-1})
proves Step~1 (a).

$(b)$
We follow the line of argument
in~\cite[Pf. of Thm.~5.1.2 Step~1]{BPS}.
Fix $a\le0$ and $\alpha\not=0$.
Then the smallest length bigger
than $a$ is given by
$\ell_\alpha>0\ge a$, where
$\ell_\alpha$ is defined
by~(\ref{eq:ell-alpha}).
Consider the \emph{new} family
$f_k=f_{k;-a,\ell_\alpha}$ as defined
subsequent to~(\ref{eq:fam-f}).
It is easy to see
that all points 
of slope $\lambda\in\Lambda_\alpha$
on the graph of $f_k$
(which is nonpositive on $\R_0^+$),
are located strictly below the line
$r\mapsto -a+\ell_\alpha r$
(which is strictly positive
on $\R^+$ due to $-a\ge0$ and
$\ell_\alpha>0$).
Therefore they are located
strictly below the line
$r\mapsto -a+\lambda r$.
Hence, by method 1
of Remark~\ref{ex:main},
all 1-periodic orbits (representing
$\alpha\not=0$)
are of action
strictly bigger than $a$.
This shows that
the action window $(a,\infty)$
is admissible and
$Spec(f_k;\alpha)\subset(a,\infty)$.
Consider the upper horizontal row
in diagram~(\ref{eq:comm-diag1})
in case
$(a,b,c)=(-\infty,a,\infty)$
and $H^0=f_k$.
It follows $[\iota^F]=0$
and
$$
     [\pi^F]:
     HF_*^{(a,\infty)}(f_k;\alpha)
     \to
     HF_*^{(-\infty,\infty)}
     (f_k;\alpha)
$$
is an isomorphism.
The target space is
independent of the Hamiltonian,
by Floer continuation,
and therefore zero: Every
sufficiently $C^2$-small
compactly supported
Hamiltonian admits only
\emph{contractible}
1-periodic orbits.
The same exhausting sequence argument
as in part $(a)$ proves~$(b)$
and therefore Step~1.
\end{proof}

\vspace{.1cm}
Fix $a,c>0$.
Step~2 and Step~3 of the proof rely
on a family of functions
$h_\delta\in
\Hh^{a,\infty}_{\alpha;c}(DT^*M,M)$
defined as follows.
Choose a real number
$\delta$ such that
\begin{equation}\label{eq:delta}
     \delta\in\left(0,\frac{a}{c}\right),\qquad
     \mu_\delta:=\frac{a}{\delta}-(c-a)
     \notin\Lambda_\alpha.
\end{equation}
It follows $\mu_\delta>0$.
Using the conventions
$\sup\emptyset=0$ and
$\inf\emptyset=\infty$,  we define
\begin{equation}\label{eq:ell}
     \ell_-:=\sup\left((0,\mu_\delta)
     \cap\Lambda_\alpha\right),\quad
     \ell_+:=\inf\left((\mu_\delta,\infty)
     \cap\Lambda_\alpha\right),\quad
     \ell_-':=\frac{\mu_\delta+\ell_-}{2}.
\end{equation}
Since $\mu_\delta\notin\Lambda_\alpha$,
Lemma~\ref{le:length-spec} shows
$$
     0\le\ell_-<\ell_-'
     <\mu_\delta<\ell_+\le+\infty,
     \qquad (\ell_-,\ell_+)\cap\Lambda_\alpha
     =\emptyset.
$$
Assumption
$\delta<a/c$ ensures that the lines
$-c+(a/\delta)r$ and $-a+\mu_\delta r$
intersect strictly above the $r$-axis
(see Figure~\ref{fig:fig-c-1}).
As a first approximation
for $h_\delta$
consider the piecewise linear
curve in $\R^2$
obtained by starting
at the point $(0,-c)$ with
slope $a/\delta$.
Upon meeting the horizontal
line through
$(0,\max\{0,\ell_-'-a\})$, follow
this horizontal line to the
right until
it intersects the vertical
line through
$(1,0)$. Now go straight down
to the point $(1,0)$
and follow the horizontal
coordinate axis
to $+\infty$. Smooth out
this piecewise
linear curve near its corners such
that the result $h_\delta$
is the graph of a smooth
function, as indicated in
Figure~\ref{fig:fig-c-1},
of slope zero
at $r=0$ and at $r=1$.
\begin{figure}[ht]
  \centering
  \epsfig{figure=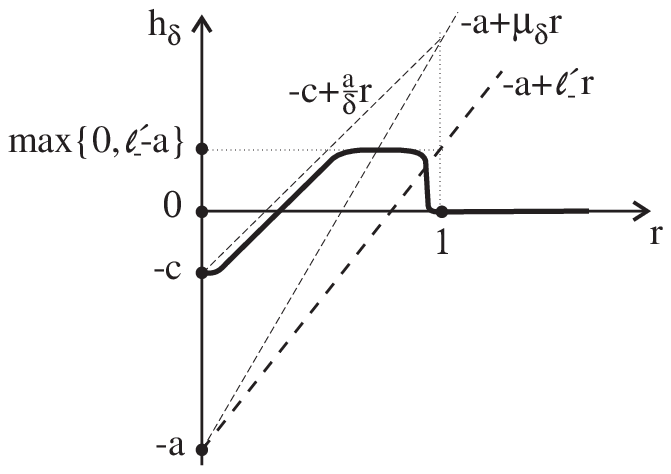}
  \caption{The function $h_\delta$
           in case $c<a$.}
  \label{fig:fig-c-1}
\end{figure}
%

\vspace{.15cm}
\noindent
{\sc Step 2.}
{\it If $a>c>0$, then
$\underrightarrow{SH}_*^{(a,\infty);c}
(DT^*M,M;\alpha)=0$.
}

\begin{proof}
Let $\delta$ be as in~(\ref{eq:delta})
and $h_\delta$ the associated function
(see Figure~\ref{fig:fig-c-1}),
then
\begin{equation}\label{eq:claim-1-S2}
     HF_*^{(a,\infty)}(h_\delta;\alpha)
     = 0.
\end{equation}
The reason is that the action
of all 1-periodic orbits
is strictly less than $a$.
This is because all
tangents to the graph
of $h_\delta$ intersect
the vertical coordinate axis
strictly above $-a$
(see Remark~\ref{ex:main} method 2).
Now choose a strictly decreasing
real sequence $\delta_k \to 0$
such that each $\delta_k$
satisfies~(\ref{eq:delta}).
Then $\{h_k:=h_{\delta_k}\}_{k\in\N}$
is an upward exhausting sequence for
the direct limit.
We have to check two properties
(see~\cite[Section 4.7]{BPS}):
Firstly, given any
$H\in\Hh^{a,\infty}_{\alpha;c}(DT^*M,M)$,
there exists $h_k\ge H$.
This is clearly true.
Secondly, the monotone homomorphism
$$
     \sigma_{h_{k+1},h_k}:
     HF_*^{(a,\infty)}(h_k;\alpha)
     \to HF_*^{(a,\infty)}(h_{k+1};\alpha),
$$
is an isomorphism, $\forall k\in\N$.
This holds trivially, since
both Floer homologies
are zero by~(\ref{eq:claim-1-S2}).
By Lemma~4.7.1~(i) in \cite{BPS},
the homomorphism in~(\ref{eq:iota})
$$
     \iota_{h_k}:
     HF_*^{(a,\infty)}(h_k;\alpha)
     \to
     \underrightarrow{SH}_*^{(a,\infty);c}
     (DT^*M,M;\alpha),\qquad
     x_{h_k}\mapsto [h_k,x_{h_k}],
$$
is an isomorphism,
$\forall k\in\N$.
By~(\ref{eq:claim-1-S2}),
the proof of Step 2 is complete.
\end{proof}

\vspace{.1cm}
\noindent
{\sc Step 3.}
{\it Given $0<a\le c$ and
$\delta$
satisfying~(\ref{eq:delta}),
there is a natural isomorphism
\begin{equation}\label{eq:step-3}
     \phi_*^{(\delta)}:
     HF_*^{(a,\infty)}(h_\delta;\alpha)
     \to
     H_*(\Ll_\alpha^{\mu_\delta^2/2} M)
     ,\qquad
     \mu_\delta
     :=\frac{a}{\delta}-(c-a).
\end{equation}
}

\begin{proof}
\begin{figure}[ht]
  \centering
  \epsfig{figure=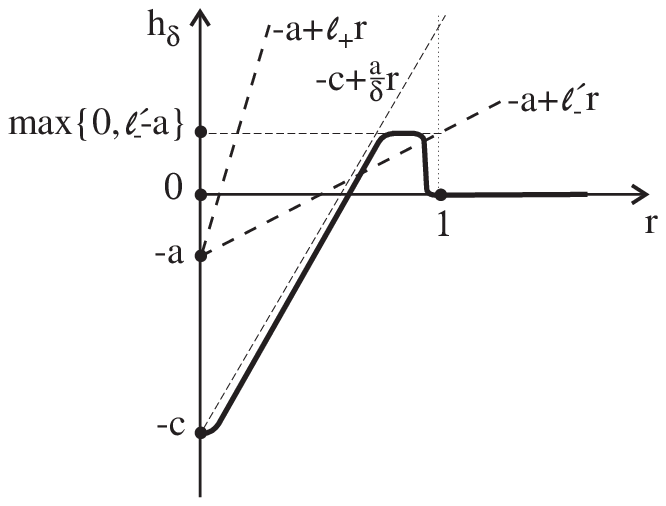}
  \caption{The function $h_\delta$ in 
           case $c\ge a$.}
  \label{fig:fig-b-1}
\end{figure}
The idea is
to deform $h_\delta$
(see Figure~\ref{fig:fig-b-1})
by action-regular
monotone homotopies
to a convex function, so that
Theorem~\ref{thm:convex-hamiltonians}
applies.
Let $0\le\ell_-<\ell_-'
<\mu_\delta<\ell_+\le+\infty$
be given by~(\ref{eq:ell}).
First we check that
the action window $(a,\infty)$
is indeed admissible for $h_\delta$:
We need to check only positive
critical slopes.
They appear in two clusters,
one of which is located near the point
$(0,-c)$. Their tangents hit the vertical
coordinate axis below $-c<-a$. So they
correspond to action values strictly
bigger than $a$. The other cluster
is located in the region between
the lines $-a+\ell_- r$ and 
$-a+\ell_+ r$.
The slope of any such tangent running
through $(0,-a)$ lies in the interval
$(\ell_-,\ell_+)$
and is therefore noncritical.
Here we use the assumption
$\mu_\delta\notin\Lambda_\alpha$,
which implies
$\ell_-<\mu_\delta<\ell_+$
and $(\ell_-,\ell_+)\cap
\Lambda_\alpha=\emptyset$.

Next we define the graph of
a new function
$\tilde{h}_\delta$ by initially
following the graph of $h_\delta$.
When $h_\delta$ makes
its first right turn,
keep going with constant
slope $a/\delta$
until you are about to hit the line
$-a+\mu_\delta r$.
Perform a smooth right turn
and follow that line
closely and linearly
as indicated in
Figure~\ref{fig:fig-b-2}.
\begin{figure}[h]
  \centering
  \epsfig{figure=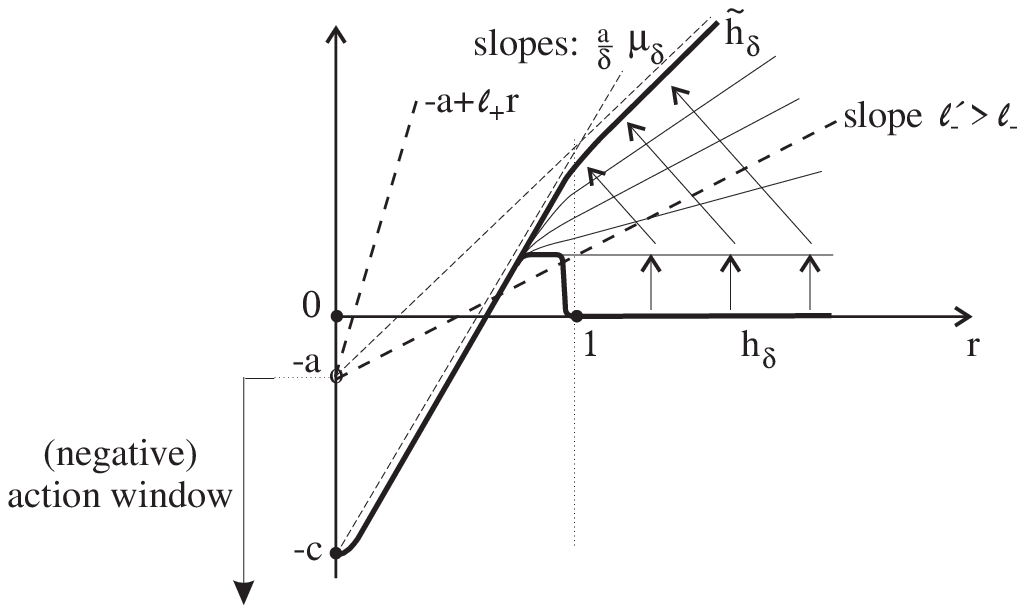}
  \caption{An action-regular monotone
           homotopy
           between $h_\delta$ and
           $\tilde{h}_\delta$.}
  \label{fig:fig-b-2}
\end{figure}
The figure also shows an
action-regular
monotone homotopy
between $h_\delta$ and
$\tilde{h}_\delta$:
We need to make sure
that no 'critical slope tangents'
to members of the homotopy move
over the point $(0,-a)$, because this
corresponds to periodic
orbits entering
or leaving the action window
$(a,\infty)$.
The former is true, since
all points on
members of the homotopy
whose tangents run through
$(0,-a)$ lie strictly
between the lines $-a+\ell_-r$
and $-a+\ell_+r$
(see Figure~\ref{fig:fig-b-2}).
Hence these slopes
are noncritical and we obtain
the monotone isomorphism
\begin{equation}\label{eq:phi1}
     \sigma_{\tilde{h}_\delta h_\delta}:
     HF_*^{(a,\infty)}(h_\delta;\alpha)
     \to
     HF_*^{(a,\infty)}
     (\tilde{h}_\delta;\alpha).
\end{equation}

Define the function
$h_\delta^{(\mu_\delta)}$
by following
$\tilde{h}_\delta$
until slope $\mu_\delta$
shows up for the first time.
Then make a smooth
turn and continue linearly with
slope $\mu_\delta$ as indicated
in Figure~\ref{fig:fig-b-3}.
\begin{figure}[h]
  \centering
  \epsfig{figure=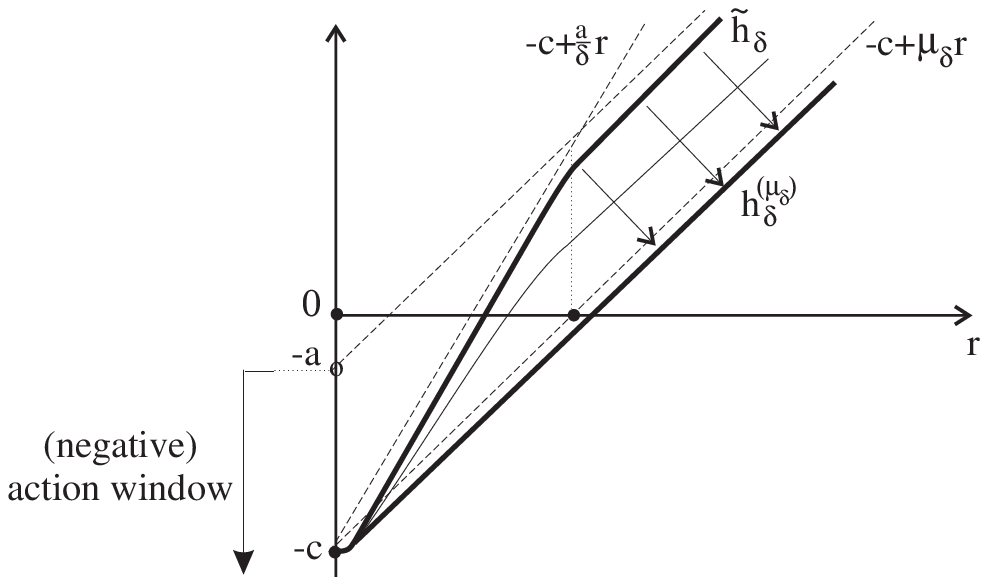}
  \caption{An action-regular monotone
           homotopy between
           $h_\delta^{(\mu_\delta)}$
           and
           $\tilde{h}_\delta$.}
  \label{fig:fig-b-3}
\end{figure}
The monotone homotopy shown
is also action-regular:
All tangents to all members of the
homotopy hit the vertical
coordinate axis
strictly below the point $-a$.
Hence the corresponding monotone
homomorphism is an 
isomorphism with inverse
\begin{equation}\label{eq:phi2}
     {\sigma_{\tilde{h}_\delta
     h_\delta^{(\mu_\delta)}}}^{-1}
     :
     HF_*^{(a,\infty)}
     (\tilde{h}_\delta;\alpha)
     \to
     HF_*^{(a,\infty)}
     (h_\delta^{(\mu_\delta)};\alpha).
\end{equation}
Since $h_\delta^{(\mu_\delta)}$
does not admit 1-periodic
orbits of action less than $a$,
the map $[\pi^F]$
in~(\ref{eq:comm-diag1})
with $(a,b,c,)=(-\infty,a,\infty)$
is an isomorphism with inverse
\begin{equation}\label{eq:phi3}
     (\pi^F)^{-1}:
     HF_*^{(a,\infty)}
     (h_\delta^{(\mu_\delta)};\alpha)
     \to
     HF_*^{(-\infty,\infty)}
     (h_\delta^{(\mu_\delta)};\alpha).
\end{equation}
Similarly, there is
the isomorphism
\begin{equation}\label{eq:phi4}
     [\pi^F]^{-1}:
     HF_*^{(-\infty,\infty)}
     (h_\delta^{(\mu_\delta)};\alpha)
     \to
     HF_*^{(-\infty,\hat{c})}
     (h_\delta^{(\mu_\delta)};\alpha),
\end{equation}
where the constant
$\hat{c}=\hat{c}(h_\delta^{(\mu_\delta)},
\mu_\delta)$ is defined in
Theorem~\ref{thm:convex-hamiltonians}.
Using again
$\mu_\delta\notin\Lambda_\alpha$,
we may apply the theorem
to obtain the isomorphism
\begin{equation}\label{eq:phi5}
     \Psi^{\mu_\delta}:
     HF_*^{(-\infty,\hat{c})}
     (h_\delta^{(\mu_\delta)};\alpha)
     \to
     H_*(\Ll_\alpha
     ^{{\mu_\delta}^2/2} M).
\end{equation}
The proof of Step~3 concludes by
composing~(\ref{eq:phi1}-\ref{eq:phi5}).
\end{proof}

\vspace{.15cm}
\noindent
{\sc Step 4.}
{\it If $0<a\le c$, then
$
     \underrightarrow{SH}_*
     ^{(a,\infty);c}
     (DT^*M,M;\alpha)\simeq
     H_*(\Ll_\alpha M)
$
naturally.
}
\begin{proof}
Consider the set
$\Delta:=\{\delta\in(0,a/c)\mid
\text{$(\ref{eq:delta})$ holds}\}$
with partial order defined
by $\delta_0\preccurlyeq\delta_1$
iff $\delta_0\ge\delta_1$.
Observe that
$(\Delta,\preccurlyeq)$
is upward directed and that
$\delta_0\ge\delta_1$
is equivalent to
$\mu_{\delta_0}\le\mu_{\delta_1}$
and to
$h_{\delta_0}\le h_{\delta_1}$.
Setting
$\mu_\Delta:=\{\mu_\delta\mid
\delta\in\Delta\}$ and
$\Hh_\Delta:=\{h_\delta\mid
\delta\in\Delta\}$,
we can therefore identify
the upward directed sets
$(\Delta,\preccurlyeq)$,
$(\mu_\Delta,\le)$
and $(\Hh_\Delta,\le)$.
Taking the direct limit
of both sides of~(\ref{eq:step-3})
with respect to
$(\Delta,\preccurlyeq)$ leads to
$$
     \underset{h_\delta\in\Hh_\Delta}
     {\underrightarrow{\lim}}
     HF_*^{(a,\infty)}(h_\delta;\alpha)
     \simeq
     \underset{\mu\in\mu_\Delta}
     {\underrightarrow{\lim}}
     H_*(\Ll_\alpha^{\mu^2/2} M).
$$
By Lemma~\ref{le:sing-hom},
the right hand side is
naturally isomorphic to 
$H_*(\Ll_\alpha M)$.
Concerning the left hand
side, there is the natural
inclusion homomorphism
$$
     \underset{h_\delta\in\Hh_\Delta}
     {\underrightarrow{\lim}}
     HF_*^{(a,\infty)}(h_\delta;\alpha)
     \to
     \underset{H\in\Hh^{a,\infty}
     _{\alpha;c}}
     {\underrightarrow{\lim}}
     HF_*^{(a,\infty)}(H;\alpha),\qquad
     [h_\delta,x_{h_\delta}]
     \mapsto
     [h_\delta,x_{h_\delta}].
$$
Using~(\ref{eq:mon-hom-comp})
and the fact that $(\Hh_\Delta,\le)$
is upward directed,
it is easy to see that
this map is well
defined, injective and surjective.
\end{proof}
Steps 2 and 4 prove part $(ii)$
of Theorem~\ref{thm:5.1.2}.

\vspace{.15cm}
\noindent
{\sc Step 5.}
{\it We prove part $(iii)$
of Theorem~\ref{thm:5.1.2}.}

\begin{proof}
Throughout we use the notation
introduced in Steps~1-4.
Fix $0<a\le c$ and choose 
$\delta\in\Delta$
sufficiently small such that
$a)$ $\mu_\delta>a$, and
$k\in\N$ sufficiently large
such that $b)$
$f_k\in\Hh^{a,\infty}_{\alpha;c}$
and $c)$ $f_k\le h_\delta$.
Let $\psi^{(k)}$ and
$\phi^{(\delta)}$ be given
by~(\ref{goal-1})
and~(\ref{eq:step-3}), respectively,
and assume that
the following diagram commutes
\begin{equation}\label{eq:comm-claim}
\begin{split}
\xymatrix{
     HF^{(a,\infty)}_*(f_k;\alpha)
     \ar[r]^{\sigma_{h_\delta f_k}}
     \ar[d]_{\psi^{(k)}}^\simeq
    &
     HF^{(a,\infty)}_*
     (h_\delta;\alpha)
     \ar[d]_\simeq^{\phi^{(\delta)}}
    \\
     H_*(\Ll^{a^2/2}_\alpha M)
     \ar[r]^{[\iota_{
     \frac{\mu_\delta^2}{2}
     \frac{a^2}{2}}]}
    &
     H_*(\Ll^{\mu_\delta^2/2}_\alpha
     M)
}
\end{split},
\end{equation}
Then the subsequent diagram,
whose existence uses~$b)$,
commutes too (thereby proving
Step~5)
\begin{equation}\label{eq:claim-iii}
\begin{split}
\xymatrix{
     \underleftarrow{SH}_*
     ^{(a,\infty)}(DT^*M;\alpha)
     \ar[rd]^{T^{(a,\infty)}_{\alpha;c}}
     \ar[d]_{\pi_{f_k}}^\simeq
    &
    \\
     HF^{(a,\infty)}_*(f_k;\alpha)
     \ar[r]^{\iota_{f_k}\qquad}
     \ar[d]_{\psi^{(k)}}^\simeq
    &
     \underrightarrow{SH}_*
     ^{(a,\infty);c}
     (DT^*M,M;\alpha)
     \ar[d]_\simeq^\phi
    \\
     H_*(\Ll^{a^2/2}_\alpha M)
     \ar[r]^{I_{a^2/2}\qquad\qquad}
    &
     \underset{\mu\in\mu_\Delta}
     {\underrightarrow{\lim}}
     H_*(\Ll_\alpha^{\mu^2/2} M)
     \simeq H_*(\Ll_\alpha M)
}
\end{split}.
\end{equation}
To see this observe that
the upper triangle
commutes by~(\ref{eq:nat-homo}).
Definition of $\phi$ and $I_{a^2/2}$
is obvious.
Commutativity of the lower
rectangular block follows
by applying~(\ref{eq:comm-claim})
to the identities
$$
     I_{a^2/2}(\psi^{(k)} x_{f_k})
     =\bigl[\mu_\delta^2/2,
     [\iota_{
     \frac{\mu_\delta^2}{2}
     \frac{a^2}{2}}]
     (\psi^{(k)} x_{f_k})\bigr],
$$
which uses $a)$, and
\begin{equation*}
\begin{split}
     \phi\circ\iota_{f_k}(x_{f_k})
     =\phi ([f_k,x_{f_k}])
     =\phi ([h_\delta,
     \sigma_{h_\delta f_k}x_{f_k}])
     =[\mu_\delta^2/2,
     \phi^{(\delta)}(
     \sigma_{h_\delta f_k}x_{f_k})],
\end{split}
\end{equation*}
which uses $c)$.
It remains to prove that
diagram~(\ref{eq:comm-claim})
commutes, which we rewrite
as follows (simplifying
notation slightly)
\begin{equation*}
\begin{split}
\xymatrix{
     HF^{(a,\infty)}_*(f_k)
     \ar[r]^\sigma
     \ar[d]_\sigma^{(\ref{eq:psi1})}
    &
     HF^{(a,\infty)}_*(h_\delta)
     \ar[d]^\sigma_{(\ref{eq:phi1})}
    \\
     HF^{(a,\infty)}_*(\tilde{f}_k)
     \ar[r]^\sigma
    &
     HF^{(a,\infty)}_*(\tilde{h}_\delta)
    \\
     HF^{(a,\infty)}_*(f_k^{(a)})
     \ar[u]^\sigma_{(\ref{eq:psi2})}
     \ar[r]^\sigma
    &
     HF^{(a,\infty)}_*(h_\delta^{(\mu_\delta)}
     )
     \ar[u]_\sigma^{(\ref{eq:phi2})}
    \\
     HF^{(-\infty,\infty)}_*(f_k^{(a)})
     \ar[u]^{[\pi^F]}_{(\ref{eq:psi3})}
     \ar[r]^\sigma
     \ar[dr]^{\Check{\sigma}}
    &
     HF^{(-\infty,\infty)}_*(h_\delta^{(\mu_\delta)}
     )
     \ar[u]_{[\pi^F]}^{(\ref{eq:phi3})}
    \\
     HF^{(-\infty,c_{f,a})}_*
     (f_k^{(a)})
     \ar[u]^{[\iota^F]}_{(\ref{eq:psi4})}
     \ar[r]^{\Check{\sigma}=\Phi^a_{hf}}_\simeq
     \ar[dr]_{\Phi^a_f,(\ref{eq:psi5})}
    &
     HF^{(-\infty,c_{h,a})}_*(
     h_\delta^{(\mu_\delta)})
     \ar[u]^{[\iota^F]}
     \ar[r]^{[\iota^F]}
     \ar[d]^{\Phi^a_h}_\simeq
    &
     HF^{(-\infty,c_{h,\mu})}_*(
     h_\delta^{(\mu_\delta)})
     \ar[ul]_{[\iota^F],(\ref{eq:phi4})}
     \ar[d]^{\Phi^{\mu_\delta}_h}_{(\ref{eq:phi5})}
    \\
    &
     H_*(\Ll^{a^2/2}_\alpha M)
     \ar[r]^{[\iota_{\mu a}]}
    &
     H_*(\Ll^{\mu_\delta^2/2}_\alpha M)
}
\end{split}.
\end{equation*}
The first two rectangular
blocks of the
diagram
commute due to
property~(\ref{eq:mon-hom-comp})
of the monotone homomorphisms.
The third block
commutes by~(\ref{eq:comm-diag1}).
The fourth block,
which consists of two triangles,
commutes
by~(\ref{eq:comm-diag2})
with $(a,b_1,b_2,c)=(-\infty,
c_{h,a},c_{f,a},\infty)$.
The triangle to its right
commutes already
on the chain level.
The final rectangular
block commutes by
Theorem~\ref{thm:convex-hamiltonians}
~(\ref{eq:thm-com2})
and the triangle to its left
by~(\ref{eq:thm-com1}).
This concludes the proof of Step~5.
\end{proof}
The proof of
Theorem~\ref{thm:5.1.2}
is complete.
\end{proof}


\section{Relative BPS-capacity}
\label{sec:BPS-capacity}

We compute the BPS-capacity
of $DT^*M$ relative to $M$
and show that every nontrivial
free homotopy class of loops in
$M$ is symplectically essential.
For $c>0$, set
$
     \Hh_c
     :=\{H\in C_0^\infty(S^1\times DT^*M)\mid
     \sup_{S^1\times M} H\le-c\}
$.
In the following definition
we use the conventions
$\inf\emptyset=\infty$ and 
$\sup\emptyset=-\infty$.
\begin{definition}
[{\cite[Sections~3.2 and~4.9]{BPS}}]
\rm
For $\alpha\in\pi_1(M)$ and $a\in\R$,
the
\emph{relative Biran-Polterovich-Salamon capacity}
is defined by
\begin{equation*}
\begin{split}
     c_{BPS}(DT^*M,M;\alpha,a)
    &:=\inf\{c>0\mid \text{$\forall
     H\in\Hh_c$ $\exists
     z\in\Pp(H;\alpha)$} \\
    &\phantom{:=\inf\{c>0\mid}
     \quad\text{such that
     $\Aa_H(z)\ge a$}\} \\
     c_{BPS}(DT^*M,M;\alpha)
    &:=\inf\{c>0\mid\text{
     $\Pp(H;\alpha)\not=\emptyset$
     $\forall H\in\Hh_c$}\}.
\end{split}
\end{equation*}
For $c>0$, define the set
$$
     A_c(DT^*M,M;\alpha)
     :=\{\text{$a\in\R$ ($a>0$
     if $\alpha=0$) $\mid$
     $T_{\alpha;c}^{(a,\infty)}
     \not=0$}\}.
$$
Define the \emph{homological
relative BPS-capacity} by
$$
     \hat{c}_{BPS}(DT^*M,M;\alpha,a)
     :=\inf\{c>0\mid
     \sup A_c(DT^*M,M;\alpha)>a \}.
$$
\end{definition}
The significance of the homological
relative capacity lies in the fact
that it is accessible to
computation and bounds
$c_{BPS}$ from
above~\cite[Proposition~4.9.1]{BPS}.
\begin{proposition}\label{pr:c-hat-BPS}
$\hat{c}_{BPS}(DT^*M,M;\alpha,a)
=\max\{\ell_\alpha,a\}$,
$\forall \alpha\in\pi_1(M)$,
$\forall a\in\R$.
\end{proposition}
\begin{proof}
Theorem~\ref{thm:5.1.2} shows that
$$
     \sup A_c(DT^*M,M;\alpha)
     =\sup [\ell_\alpha,c]
     =
     \begin{cases}
       c,
      &\text{if $c\ge\ell_\alpha$},\\
       -\infty,
      &\text{else,}
     \end{cases}
$$
and this implies
$$
     \hat{c}_{BPS}(DT^*M,M;\alpha,a)
     =\inf\{c>0\mid
     \text{$c>a$ and $c\ge\ell_\alpha$}\}
     =\max\{a,\ell_\alpha\}.
$$
\end{proof}
\begin{theorem}[BPS-capacity of
\boldmath$DT^*M$ relative $M$]
\label{thm:3.2.1}
Let $M$ be a closed connected smooth Riemannian
manifold. Then, for every
free homotopy class $\alpha$ of loops in $M$
and every $a\in\R$, the
relative BPS-capacity
is finite and given by
$$
     c_{BPS}(DT^*M,M;\alpha,a)
     =\max\{\ell_\alpha,a\}.
$$
\end{theorem}
\begin{proof}
The proof
of~\cite[Proof of Theorem~3.2.1]{BPS}
carries over almost literally.
In case $\alpha=0$ every
$H\in\Hh$ admits constant
1-periodic orbits of action zero,
because it is zero near the
boundary of $DT^*M$.
Therefore $c_{BPS}(DT^*M,M;0,a)=0$
whenever $a\le0$.
The remaining part of proof proceeds
like~\cite[Proof of Theorem~3.2.1; case $\alpha\not=0$]{BPS}
with Theorem~5.1.1 and
Theorem~5.1.2 replaced by
Proposition~\ref{pr:c-hat-BPS},
$\ell$ by $\ell_\alpha$,
and with
$$
     m
     :=
     \begin{cases}
       \max\{\ell_\alpha,a\},
      &\text{in case $\alpha\not=0$ and $a\in\R$,}\\
       a,
      &\text{in case $\alpha=0$ and $a>0$.}
     \end{cases}
$$
\end{proof}

\begin{definition}
[{\cite[Section~3.4]{BPS}}]
\label{def:symp-ess}\rm
Let $M$ be a closed connected smooth manifold.
A nontrivial free homotopy class
$\alpha$ of loops in $M$ is called
\emph{symplectically essential}
if there exists a domain $W\subset T^*M$
containing $M$ and such that
$c_{BPS}(W,M;\iota_\#\alpha)$ is finite.
Here $\iota_\#$ is the map between
free homotopy classes induced by the inclusion
$\iota:M\hookrightarrow W$.
\end{definition}

\begin{corollary}\label{cor:symp-ess}
Every nontrivial free homotopy class $\alpha$
of loops in a closed connected smooth Riemannian manifold
$M$ is symplectically essential.
\end{corollary}

\begin{proof}
Let $W:=DT^*M$, then $\iota_\#$
is an isomorphism and in abuse of
notation we write $\iota_\#\alpha=\alpha$.
Theorem~\ref{thm:3.2.1} shows
\begin{equation}\label{eq:cBPS}
     c_{BPS}(DT^*M,M;\alpha)
     =c_{BPS}(DT^*M,M;\alpha,-\infty)
     =\ell_\alpha.
\end{equation}
\end{proof}

\begin{proof}[Proof of Theorem~\ref{thm:exist-orbit}]
By~\cite[Proof of Theorem~B]{BPS},
we may assume without loss of
generality that $H$ is periodic, i.e.
$H\in C^\infty_0 (S^1\times DT^*M)$.
Theorem~\ref{thm:3.2.1} shows
$
     c_{BPS}(DT^*M,M;\alpha,c)
     =\max\{\ell_\alpha,c\}
     =c
$.
By~\cite[Proposition~3.3.4]{BPS},
the set which appears in the definition
of $c_{BPS}(DT^*M,M;\alpha,c)$
has a minimum
and this proves the theorem.
\end{proof}






\begin{thebibliography}{99}
\small

\bibitem{Bo82}  R.~Bott,
      Lectures on Morse theory, old and new,
      {\it Bull. (New series) Amer. Math. Soc.}
      {\bf 7} (1982), 331--358.

\bibitem{BPS}  P.~Biran, L.~Polterovich
      and D.A.~Salamon,
      Propagation in Hamiltonian dynamics
      and relative symplectic homology,
      {\it Duke Math. Journal}
      {\bf 119} (2003), 65--118.

\bibitem{Ci94}  K.~Cieliebak,
      Pseudo-holomorphic curves and
      periodic orbits on cotangent bundles,
      {\it J. Math. Pures Appl.} {\bf 73}
      (1994), 251--278.

\bibitem{CFH}  K.~Cieliebak, A.~Floer and H.~Hofer,
      Symplectic homology,
      II. A general construction,
      {\it Math. Zeit.} {\bf 218} (1995), 103--122.

\bibitem{FLOER5}  A.~Floer,
      Symplectic fixed points and
      holomorphic spheres,
      {\it Comm. Math. Phys.} {\bf 120} (1989),
      575--611.

\bibitem{FH1}  A.~Floer and H.~Hofer,
      Symplectic homology,
      I. Open sets in $\C^n$,
      {\it Math. Zeit.} {\bf 215} (1994), 37--88.

\bibitem{FHS}  A.~Floer, H.~Hofer and
      D.A.~Salamon,
      Transversality in elliptic Morse theory
      for the symplectic action,
      {\it Duke Math. Journal}
      {\bf 80} (1995), 251--292.

\bibitem{Gi03}  V.L.~Ginzburg,
      The Weinstein conjecture and the
      theorems of nearby and almost existence,
      Preprint 2003, math.DG/0310330.

\bibitem{GL}  D.~Gatien and F.~Lalonde,
      Holomorphic cylinders with
      Lagrangian boundaries and
      Hamiltonian dynamics,
      {\it Duke Math. Journal}
      {\bf 102} (2000), 485--511.

\bibitem{GT77}  D.~Gilbarg and N.S.~Trudinger,
      {\it Elliptic partial differential equations 
      of second order},
      Grundlehren der mathematischen
      Wissenschaften {\bf 224}, 
      Springer-Verlag 1977, third printing 1998.

\bibitem{HV88}  H.~Hofer and C.~Viterbo,
      The Weinstein conjecture in cotangent
      bundles and related results,
      {\it Annali Sc. Norm. Sup. Pisa},
      Serie IV, Fasc. III,
      {\bf 15} (1988), 411--445.

\bibitem{HZ94}  H.~Hofer and E.~Zehnder,
      {\it Symplectic invariants
      and Hamiltonian dynamics},
      Birkh\"auser Advanced Texts,
      Basel, 1994.

\bibitem{JO}  J.~Jost,
      {\it Riemannian geometry and
      geometric analysis},
      Universitext, Springer Verlag,
      Berlin Heidelberg, 1995.

\bibitem{Lee}  Y.-J.~Lee,
      Non-contractible periodic orbits,
      Gromov invariants, and
      Floer-theoretic torsions,
      Preprint 2003, math.SG/0308185.

\bibitem{Oa03}  A.~Oancea,
      PhD thesis, Universite Paris XI,
      UFR Scientifique d'Orsay, 2003. 

\bibitem{Sa97}  D.A.~Salamon,
      Lectures on Floer Homology,
      In {\it Symplectic Geometry and Topology},
      edited by Y.~Eliashberg and L.~Traynor,
      IAS/Park City Mathematics Series,
      Vol {\bf 7}, 1999, pp 143--230.

\bibitem{SW1}  D.A.~Salamon and J.~Weber,
      Floer homology and the heat flow,
      Preprint, ETH-Z\"urich, April 2003.

\bibitem{SALZ}  D.A.~Salamon and E.~Zehnder,
      Morse theory for periodic solutions of Hamiltonian systems and
      the Maslov index, {\it Comm. Pure Appl. Math.}
      {\bf 45} (1992), 1303--1360.

\bibitem{Sa42}  A.~Sard,
      The measure of the critical points
      of differentiable maps,
      {\it Bull. Amer. Math. Soc.}
      {\bf 48} (1942), 883--890.

\bibitem{Vi87}  C.~Viterbo,
      A proof of Weinstein's conjecture
      in $\R^{2n}$,
      {\it Ann. Inst. Poincar\'e,
      Anal. Non Lin\'eaire},
      {\bf 4} (1987), 337--356.

\bibitem{V}  C.~Viterbo,
      Functors and computations in
      Floer homology with
      applications, I,
      {\it Geom. funct. anal.}
      {\bf 9} (1999), 985--1033.

\bibitem{Vi97}  C.~Viterbo,
      Exact Lagrange submanifolds,
      periodic orbits and the
      cohomology of free loop spaces,
      {\it J. Diff. Geom.}
      {\bf 47} (1997), 420--468.

\bibitem{Vi98}  C.~Viterbo,
      Functors and computations in Floer homology
      with applications, Part II,
      Preprint~1996, revised~1998.

\bibitem{JOA1}  J.~Weber,
      Perturbed closed geodesics are periodic orbits:
      Index and transversality, 
      {\it Math. Z.} {\bf 241}
      (2002), 45--81. 
      http://dx.doi.org/10.1007/s002090100406.

\bibitem{We79}  A.~Weinstein,
      On the hypothesis of Rabinowitz'
      periodic orbit theorems,
      {\it J. Differential Equations}
      {\bf 33} (1979), 353--358.

\end{thebibliography}
\end{document}